\documentclass{article}
\usepackage{pb-diagram}
\usepackage{pst-grad}
\usepackage{pstricks,}
\usepackage{pst-plot}
\usepackage[tiling]{pst-fill}
\usepackage{amsfonts, amssymb, amsmath}
\usepackage[all]{xy}

\usepackage[bookmarks=false,breaklinks=true]{hyperref}
\usepackage{color}
\definecolor{myurlcolor}{rgb}{0,0,0.8}
\hypersetup{colorlinks,
linkcolor=myurlcolor,
citecolor=myurlcolor,
urlcolor=myurlcolor}
\usepackage[mathscr]{euscript}

\def\stackto #1 { \, {\stackrel{#1}{\longrightarrow}}\, }
\def\stackTo #1 { \, {\stackrel{#1}{\Longrightarrow}}\, }

\newcommand{\et}{\hspace{-0.08in}{\bf .}\hspace{0.1in}}

\newcommand{\maps}{\colon}
\renewcommand{\to}{\rightarrow}
\newcommand{\To}{\Rightarrow}

\newcommand{\iso}{\cong}

\newcommand{\Ob}{{\rm Ob}}
\newcommand{\Mor}{{\rm Mor}}

\newcommand{\id}{{\rm id}}
\newcommand{\tr}{{\rm tr}}
\newcommand{\ad}{{\rm ad}}
\newcommand{\inv}{{\rm inv}}

\newcommand{\U}{{\rm U}}

\newcommand{\Spin}{{\rm Spin}}
\newcommand{\String}{{\rm String}}

\newcommand{\g}{{\mathfrak g}}

\renewcommand{\u}{{\mathfrak u}}
\newcommand{\lietwoalg}{{\mathfrak c}}
\newcommand{\lietwogrp}{C}

\newcommand{\LG}{{\cal L}_kG}

\newcommand{\PG}{{\cal P}_kG}
\newcommand{\Pg}{{\cal P}_k\g}
\newcommand{\OG}{\Omega G}
\newcommand{\Og}{\Omega \mathfrak{g}}
\newcommand{\wOG}{\widehat{\Omega G}}
\newcommand{\wOg}{\widehat{\Omega \mathfrak{g}}}
\newcommand{\wOkG}{\widehat{\Omega_k G}}
\newcommand{\wOkg}{\widehat{\Omega_k \g}}
\newcommand{\cG}{\mathcal{G}}

\renewcommand{\a}{\alpha}
\renewcommand{\b}{\beta}
\renewcommand{\c}{\gamma}

\newcommand{\T}{{\cal E}}

\newcommand{\R}{{\mathbb R}}

\newcommand{\Z}{{\mathbb Z}}
\newcommand{\TT}{{\U(1)}}

\newtheorem{thm}{Theorem}[section]
\newtheorem{cor}[thm]{Corollary}
\newtheorem{lem}[thm]{Lemma}

\newtheorem{prop}[thm]{Proposition}
\newtheorem{defn}[thm]{Definition}
\newtheorem{example}[thm]{Example}

\def\bracket #1#2 {{\langle {#1} , {#2} \rangle}}

\def\Proof {\noindent{\bf Proof. }}

\def\endofproof { \hskip 2em $\Box$ \vskip 1em}

\def\refdef #1 {{ (def. \ref{#1}) }}

\def\of #1 {{\!\left({#1}\right)}}

\title{From Loop Groups to 2-Groups}
     \author{John C.\ Baez and Danny Stevenson \\
     Department of Mathematics,  University of California\\
     Riverside, California 92521, USA \\
     \\
     Alissa S.\ Crans \\
     Department of Mathematics,
     Loyola Marymount University \\
     Los Angeles, CA 90045, USA \\
     \\
     Urs Schreiber \\
     Fachbereich Physik, Universit{\"a}t Duisburg-Essen \\
     Essen 45117, Germany
     \\
     \\ email: baez@math.ucr.edu, dstevens@math.ucr.edu, \\
     acrans@lmu.edu, urs.schreiber@uni-essen.de \\
     \\
     }

      \date{January 20, 2007}
\begin{document}
\bibliographystyle{plain}
\maketitle

\begin{abstract}
\noindent
We describe an interesting relation between Lie 2-algebras, the
Kac--Moody central extensions of loop groups, and the group
$\String(n)$.  A Lie 2-algebra is a categorified version of a Lie
algebra where the Jacobi identity holds up to a natural isomorphism
called the `Jacobiator'.  Similarly, a Lie 2-group is a categorified
version of a Lie group.  If $G$ is a simply-connected compact simple
Lie group, there is a 1-parameter family of Lie 2-algebras $\g_k$ each
having $\g$ as its Lie algebra of objects, but with a Jacobiator built
from the canonical 3-form on $G$.  There appears to be no Lie 2-group
having $\g_k$ as its Lie 2-algebra, except when $k = 0$.  Here,
however, we construct for integral $k$ an infinite-dimensional Lie
2-group $\PG$ whose Lie 2-algebra is {\it equivalent} to $\g_k$.  The
objects of $\PG$ are based paths in $G$, while the automorphisms of
any object form the level-$k$ Kac--Moody central extension of the loop
group $\OG$.  This 2-group is closely related to the $k$th power of
the canonical gerbe over $G$.  Its nerve gives a topological group
$|\PG|$ that is an extension of $G$ by $K(\Z,2)$.  When $k = \pm 1$,
$|\PG|$ can also be obtained by killing the third homotopy group of
$G$.  Thus, when $G = \Spin(n)$, $|\PG|$ is none other than $\String(n)$.
\end{abstract}

\section{Introduction} \label{introductionsection}

The theory of simple Lie groups and Lie algebras has long played a
central role in mathematics.  Starting in the 1980s, a wave of
research motivated by physics has expanded this theory to
include structures such as quantum groups, affine Lie
algebras, and central extensions of loop groups.  All these
structures rely for their existence on the left-invariant closed
3-form $\nu$ naturally possessed by any compact simple Lie group
$G$:
\[              \nu(x,y,z) = \langle x, [y,z] \rangle   \qquad
x,y,z \in \g , \]
or its close relative, the left-invariant closed 2-form $\omega$ on the loop
group $\OG$:
\[      \omega(f, g) = 2
\int_{S^1} \langle f(\theta), g'(\theta) \rangle \, d \theta
\qquad f,g \in \Og .\] Moreover, all these new structures fit
together in a framework that can best be understood with ideas
from physics --- in particular, the Wess--Zumino--Witten model and
Chern--Simons theory.  Since these ideas arose from work on string
theory, which replaces point particles by higher-dimensional
extended objects, it is not surprising that their study uses
concepts from higher-dimensional algebra, such as gerbes
\cite{Brylinski,BM,CJMSW}.

More recently, work on higher-dimensional algebra has focused attention
on Lie 2-groups \cite{HDA5} and Lie 2-algebras \cite{HDA6}.  A
`2-group' is a category equipped with operations analogous to those of
a group, where all the usual group axioms hold only up to specified
natural isomorphisms satisfying certain coherence laws of their own.
A `Lie 2-group' is a 2-group where the set of objects and the set of
morphisms are smooth manifolds, and all the operations and natural
isomorphisms are smooth.  Similarly, a `Lie 2-algebra' is a category
equipped with operations analogous to those of a Lie algebra, satisfying
the usual laws up to coherent natural isomorphisms.  Just as Lie
groups and Lie algebras are important in gauge theory, Lie 2-groups
and Lie 2-algebras are important in `higher gauge theory', which
describes the parallel transport of higher-dimensional extended objects
\cite{BS,Bartels}.

The question naturally arises whether every finite-dimensional Lie
2-algebra comes from a Lie 2-group.  The answer is surprisingly
subtle, as illustrated by a class of Lie 2-algebras coming from
simple Lie algebras. Suppose $G$ is a simply-connected compact
simple Lie group $G$, and let $\g$ be its Lie algebra.   For any
real number $k$, there is a Lie 2-algebra $\g_k$ for which the
space of objects is $\g$, the space of endomorphisms of any object
is $\R$, and $[[x,y], z] = [x,[y,z]] + [[x,z],y]$, but the
`Jacobiator' isomorphism
\[     J_{x,y,z} \maps [[x,y], z] \stackto{\sim} [x, [y,z]] + [[x,z],y]  \]
is not the identity; instead we have
\[      J_{x,y,z} = k \, \nu(x,y,z)  \]
where $\nu$ is as above.   If we normalize the invariant inner
product $\langle \cdot, \cdot \rangle$ on $\g$ so that the de Rham
cohomology class of the closed form $\nu/2\pi$ generates the third integral
cohomology group of $G$, then there is a 2-group $G_k$ corresponding
to $\g_k$ in a certain sense whenever $k$ is an integer \cite{HDA5}.
The construction of this 2-group is very interesting, because it uses
Chern--Simons theory in an essential way.  However, for $k \ne 0$
there is no good way to make this 2-group into a Lie 2-group!  The
set of objects is naturally a smooth manifold, and so is the set
of morphisms, and the group operations are smooth, but the
associator
\[           a_{x,y,z} \maps (xy)z \stackto{\sim} x(yz)   \]
cannot be made everywhere smooth, or even continuous.

It would be disappointing if such a fundamental Lie 2-algebra as
$\g_k$ failed to come from a Lie 2-group even when $k$ was an
integer.  Here we resolve this dilemma by finding a Lie 2-algebra
{\it equivalent} to $\g_k$ that {\it does} come from a Lie 2-group
--- albeit an infinite-dimensional one:

\begin{thm} \et \label{punchline} Let $G$ be a simply-connected compact
simple Lie group.  For any $k \in \Z$, there is a Fr\'echet Lie
2-group $\PG$ whose Lie 2-algebra $\Pg$ is equivalent to $\g_k$.
\end{thm}

\noindent
Here two Lie 2-algebras are `equivalent' if there are
homomorphisms going back and forth between them that are inverses
{\it up to natural isomorphism}.  

We also study the relation between $\PG$ and the topological group
$\hat G$ obtained by killing the third homotopy group of $G$.
When $G = \Spin(n)$, this topological group is famous under the
name of $\String(n)$, since it plays a role in string theory
\cite{MuSt,Stolz-Teichner,Witten}.  More generally, any compact
simple Lie group $G$ has $\pi_3(G) = \Z$, but after killing
$\pi_1(G)$ by passing to the universal cover of $G$, one can then
kill $\pi_3(G)$ by passing to $\hat G$, which is defined as the
homotopy fiber of the canonical map from $G$ to the Eilenberg--Mac
Lane space $K(\Z,3)$. This specifies $\hat G$ up to homotopy, but
there is still the interesting problem of finding nice geometrical
models for $\hat G$.

Stolz and Teichner \cite{Stolz-Teichner} have already given one
solution to this problem.  Here we present another.  Given any
topological 2-group $C$, the geometric realization of its nerve 
is a topological group $|C|$.  Applying this process to $\PG$ when 
$k = \pm 1$, we obtain $\hat G$:

\begin{thm} \et  Let $G$ be a simply-connected compact
simple Lie group.  Then $|\PG|$ is an extension of $G$ by a
topological group that is homotopy equivalent to $K(\Z,2)$.
Moreover, $|\PG| \simeq \hat{G}$ when $k = \pm 1$.
\end{thm}

\noindent While this construction of $\hat G$ uses simplicial
methods and is thus arguably less `geometric' than that of Stolz
and Teichner, it avoids their use of type III$_1$ von Neumann
algebras, and has a simple relation to the Kac--Moody central
extension of $G$.

The 2-group $\PG$ is easy to describe, in part because it is
`strict': all the usual group axioms hold as equations.   The
basic idea is easiest to understand using some geometry.  Apart
from some technical fine print, an object of $\PG$ is just a path
in $G$ starting at the identity.  A morphism from the path $f_1$
to the path $f_2$ is an equivalence class of pairs $(D,z)$
consisting of a disk $D$ going from $f_1$ to $f_2$:

\vskip 1em
\[
 \xy
  (0,0)*{
    \begin{pspicture}(3,3)
    \pscustom[fillcolor=lightgray,fillstyle=gradient,
      gradbegin=gray, gradend=white, gradmidpoint=1,gradangle=110]{
        \psarc(1.5,1.5){1.24}{180}{0}
        \psellipse[linestyle=none](1.5,1.5)(1.25,.5)
    }
    \pscustom[linestyle=dotted,fillcolor=lightgray,fillstyle=gradient,
            gradbegin=lightgray, gradend=white,
          gradmidpoint=1,gradangle=110]{
        \psellipse[linestyle=dotted](1.5,1.5)(1.25,.5)
        \psarc[linestyle=dotted](1.5,1.5){1.24}{0}{180}
    }
    \begin{psclip}{\pswedge[linestyle=none](1.5,1.5){1.25}{180}{0}}
          \psellipse(1.5,1.5)(1.25,.5)
    \end{psclip}
    \pscircle(1.5,1.5){1.25}
    \psdots(1.5,2.75)(2,1.2)
    \psbezier[arrows=->](1.5,2.75)(1.4,2.4)(2.4,2.3)(2,1.7)
    \psbezier(2,1.7)(1.9,1.5)(2.1,1.4)(2,1.2)
    \psbezier[arrows=->](1.5,2.75)(1,2.4)(1,1.9)(1.1,1.7)
    \psbezier(1.1,1.7)(1.3,1.3)(1.7,1.3)(2,1.2)
     \rput(-.3,2){$G$}
    \rput(1.5,2.95){$\scriptstyle 1$}
    \rput(.85,1.7){$\scriptstyle f_1$}
    \rput(2.2,1.7){$\scriptstyle f_2$}
   \end{pspicture}};
        {\ar@{=>}_<<{\scriptstyle D} (-1,3);(3.5,3) };
 \endxy
\]
\vskip 1em

\noindent together with a unit complex number $z$.  Given two such
pairs $(D_1,z_1)$ and $(D_2,z_2)$, we can always find a 3-ball $B$
whose boundary is $D_1 \cup D_2$, and we say the pairs are
equivalent when
\[   z_2/z_1 = e^{ik \int_B \nu} \]
where $\nu$ is the left-invariant closed 3-form on $G$ given as
above.  Note that $\exp(ik \int_B \nu)$ is independent of the choice
of $B$, because the integral of $\nu$ over any 3-sphere is $2\pi$
times an integer.  There is an obvious way to compose morphisms in $\PG$,
and the resulting category inherits a Lie 2-group structure from the Lie
group structure of $G$.

The above description of $\PG$ is modeled after Murray's construction
\cite{Murray} of a gerbe from an integral closed 3-form on a manifold
with a chosen basepoint.  Indeed, $\PG$ is just another way of talking
about the $k$th power of the canonical gerbe on $G$, and the 2-group
structure on $\PG$ is a reflection of the fact that this gerbe is
`multiplicative' in the sense of Brylinski \cite{Brylinski2}.  The
3-form $k\nu$, which plays the role of the Jacobiator in $\g_k$, is
the 3-curvature of a connection on this gerbe.

In most of this paper we take a slightly
different viewpoint.  Let $P_0 G$ be the space of smooth
paths $f \maps [0,2\pi] \to G$ that start at the identity of $G$.
This becomes an infinite-dimensional Lie group under pointwise
multiplication.  The map $f \mapsto f(2\pi)$ is a homomorphism from
$P_0 G$ to $G$ whose kernel is precisely $\Omega G$.  For any $k \in \Z$,
the loop group $\OG$ has a central extension
\[  1 \stackto{\;} \U(1) \stackto{\;} \wOkG \stackto{p} \OG \stackto{\;} 1 \]
which at the Lie algebra level is determined by the 2-cocycle
$ik\omega$, with $\omega$ defined as above.  This is called the
`level-$k$ Kac--Moody central extension' of $G$.
The infinite-dimensional Lie 2-group $\PG$ has $P_0 G$ as its group of
objects, and given $f_1, f_2 \in P_0 G$, a morphism $\hat{\ell}
\maps f_1 \to f_2$ is an element $\hat{\ell} \in \wOkG$ such that
\[  f_2 / f_1 = p(\hat{\ell}) . \]
In this description, composition of morphisms in $\PG$ is multiplication
in $\wOkG$, while again $\PG$ becomes a Lie 2-group using the Lie group
structure of $G$.

After we wrote the first version of this paper \cite{BCSS}, Andr\'e Henriques 
\cite{Henriques} showed that $\PG$ arises from a general theory for 
integrating Lie $n$-algebras to obtain Lie 
$n$-groups of his sort.  The two papers should be read together.

\section{Review of Lie 2-Algebras and Lie 2-Groups}

We begin with a review of Lie 2-algebras and Lie 2-groups.
More details can be found in our papers HDA5 \cite{HDA5} and
HDA6 \cite{HDA6}.  Our notation largely follows that of these papers,
but the reader should be warned that here we denote the composite of
morphisms $f \maps x \rightarrow y$ and $g \maps y \rightarrow z$ as
$g \circ f \maps x \rightarrow z.$

\subsection{Lie 2-algebras}

The concept of `Lie 2-algebra' blends together the notion of a
Lie algebra with that of a category.  Just as a Lie algebra has an
underlying vector space, a Lie 2-algebra has an underlying
2-vector space: that is, a category where everything is {\it
linear}.  

More precisely, a {\bf 2-vector space} $L$ is a category for which
the set of objects $\Ob(L)$ and the set of morphisms $\Mor(L)$
are both vector spaces, and
the maps $s, t \maps \Mor(L) \to
\Ob(L)$ sending any morphism to its source and target,
the map $i \maps \Ob(L) \to \Mor(L)$ sending any object
to its identity morphism,
and the map $\circ$ sending any composable pair of morphisms
to its composite
are all linear.  As usual, we write a
morphism as $f \maps x \to y$ when $s(f) = x$ and $t(f) = y$, and
we often write $i(x)$ as $1_x$.

To obtain a Lie $2$-algebra, we begin with a $2$-vector space
and equip it with a bracket functor, which satisfies the
Jacobi identity up to a natural isomorphism called the
`Jacobiator'. Then we require that the Jacobiator satisfy a new
coherence law of its own: the `Jacobiator identity'.

\begin{defn} \et \label{defnlie2alg}
A {\bf Lie $2$-algebra} consists of:

\begin{itemize}
\item a $2$-vector space $L$
\end{itemize}
equipped with:
\begin{itemize}
\item a functor called the {\bf bracket}
\[   [\cdot, \cdot] \maps L \times L \to L ,\]
bilinear and skew-symmetric as a function of objects and morphisms,
\item a natural isomorphism called the {\bf Jacobiator},
\[ J_{x,y,z} \maps [[x,y],z] \to [x,[y,z]] + [[x,z],y],\]
trilinear and antisymmetric as a function of the objects $x,y,z \in L$,
\end{itemize}
such that:
\begin{itemize}
\item the {\bf Jacobiator identity} holds: the following
diagram commutes for all objects $w,x,y,z \in L$:
$$ \def\objectstyle{\scriptstyle}
    \def\labelstyle{\scriptstyle}
     \xy
     (0,35)*+{[[[w,x],y],z]}="1";
     (-40,20)*+{[[[w,y],x],z] + [[w,[x,y]],z]}="2";
     (40,20)*+{[[[w,x],z],y] + [[w,x],[y,z]]}="3";
     (-40,0)*+{[[[w,y],z],x] + [[w,y],[x,z]]}="4'";
     (-40,-4)*+{+ [w,[[x,y],z]] + [[w,z],[x,y]]}="4";
     (40,0)*+{[[w,[x,z]],y]}="5'";
     (40,-4)*+{+ [[w,x],[y,z]] + [[[w,z],x],y]}="5";
     (-35,-23)*+{[[[w,z],y],x] + [[w,[y,z]],x]}="6'";
     (-35,-28)*+{+ [[w,y], [x,z]] + [w,[[x,y],z]] + [[w,z],[x,y]]}="6";
     (35,-23)*+{[[[w,z],y],x] + [[w,z],[x,y]]  + [[w,y],[x,z]]}="7'";
     (35,-28)*+{+ [w,[[x,z],y]]  + [[w,[y,z]],x] + [w,[x,[y,z]]]}="7";
          {\ar_{[J_{w,x,y},z]}                   "1";"2"};
          {\ar^{J_{[w,x],y,z}}                               "1";"3"};
          {\ar_{J_{[w,y],x,z} + J_{w,[x,y],z}}   "2";"4'"};
          {\ar_{[J_{w,y,z},x]+1}                 "4";"6'"};
          {\ar^{[J_{w,x,z},y]+1}                   "3";"5'"};
          {\ar^{J_{w,[x,z],y}+  J_{[w,z],x,y} + J_{w,x,[y,z]}}                 "5";"7'"};
          {\ar_{[w,J_{x,y,z}]+1 \; \; }          "6";"7"};
\endxy
\\ \\
$$
\end{itemize}
\end{defn}

\noindent Note that for any object $z \in L$, bracketing with $z$
defines a functor from $L$ to itself, which we use above to define
morphisms such as $[J_{w,x,y}, z]$.

A homomorphism between Lie $2$-algebras is a linear functor
preserving the bracket, but only up to a specified natural
isomorphism satisfying a suitable coherence law.  More precisely:

\begin{defn} \et \label{lie2algfunct} Given
Lie $2$-algebras $L$ and $L'$, a {\bf homomorphism} $F \maps L
\rightarrow L'$ consists of:

\begin{itemize}
    \item a functor $F$ from the underlying $2$-vector space of
        $L$ to that of $L'$, linear on objects and morphisms,
\item a natural isomorphism
        $$F_{2}(x,y)\maps [F(x), F(y)] \rightarrow F[x,y],$$
bilinear and skew-symmetric as a function of the objects $x, y \in L$,
\end{itemize}
such that:
\begin{itemize}
\item
the following diagram commutes for all objects $x,y,z \in L$:
$$\xymatrix{
      [[F(x), F(y)], F(z)]
        \ar[rrrr]^<<<<<<<<<<<<<<<<<<<<<<{J_{F(x), F(y),
F(z)}}        \ar[dd]_{[F_{2},1]}
         &&&& [F(x), [F(y), F(z)]] + [[F(x), F(z)], F(y)]
        \ar[dd]^{[1,F_{2}] + [F_{2},1]} \\ \\
         [F[x,y], F(z)]
        \ar[dd]_{F_{2}}
         &&&& [F[x], F[y,z]] + [F[x,z], F(y)]
        \ar[dd]^{F_{2} + F_{2}} \\ \\
         F[[x,y],z]]
        \ar[rrrr]^{F(J_{x,y,z})}
         &&&& F[x,[y,z]] + [F[x,z],y]}$$
\end{itemize}
\end{defn}

\noindent Here and elsewhere we omit the arguments of natural
transformations such as $F_2$ and $G_2$ when these are obvious
from context.

Similarly, a `2-homomorphism' is a linear natural isomorphism that is
compatible with the bracket structure:

\begin{defn} \et \label{lie2algnattrans} Let $F,G \maps
L \to L'$ be Lie 2-algebra homomorphisms.  A {\bf 2-homomorphism}
$\theta \maps F \To G$ is a natural transformation
\[          \theta_x \maps F(x) \to G(x) , \]
linear as a function of the object $x \in L$, such that the following
diagram commutes for all $x, y \in L$:
$$\xymatrix{
      [F(x), F(y)]
       \ar[rr]^{F_{2}}
       \ar[dd]_{[\theta_{x}, \theta_{y}]}
        && F[x,y]
       \ar[dd]^{\theta_{[x,y]}} \\ \\
        [G(x), G(y)]
       \ar[rr]^{G_{2}}
        && G[x,y] }$$

\end{defn}

In HDA6 we showed:

\begin{prop} \et There is a strict 2-category {\bf Lie2Alg}
with Lie $2$-algebras as \break objects,
homomorphisms between these as morphisms, and $2$-homomorphisms
between those as 2-morphisms.
\end{prop}

\subsection{$L_\infty$-algebras}
\label{Linfty.section}

Just as the concept of Lie 2-algebra blends the notions of Lie
algebra and category, the concept of `$L_\infty$-algebra' blends
the notions of Lie algebra and chain complex.  More precisely, an
$L_\infty$-algebra is a chain complex equipped with a bilinear
skew-symmetric bracket operation that satisfies the Jacobi
identity up to a chain homotopy, which in turn satisfies a
law of its own up to chain homotopy, and so on {\it ad infinitum}.
In fact, $L_\infty$-algebras were defined long before Lie 2-algebras,
going back to a 1985 paper by Schlessinger and Stasheff \cite{SS}.
They are also called `strongly homotopy Lie algebras',
or `sh Lie algebras' for short.

Our conventions regarding $L_\infty$-algebras follow those of Lada
and Markl \cite{LM}.  In particular, for graded objects $x_{1},
\ldots, x_{n}$ and a permutation $\sigma \in S_{n}$ we define the
{\bf Koszul sign} $\epsilon(\sigma; x_1, \ldots, x_n)$ by the
equation
\[
x_{1} \wedge \cdots \wedge x_{n} = \epsilon(\sigma; x_1, \ldots,
x_n) \, x_{\sigma(1)} \wedge \cdots \wedge x_{\sigma(n)},
\]
which must be satisfied in the free
graded-commutative algebra on $x_{1}, \ldots, x_{n}.$
Furthermore, we define
\[
\chi(\sigma) = \textrm{sgn} (\sigma) \,
\epsilon(\sigma; x_{1}, \dots, x_{n}).
\]
Thus, $\chi(\sigma)$ takes into account the sign of the permutation in
$S_{n}$ as well as the Koszul sign.
Finally, if $n$ is a natural number and $1 \leq j \leq n-1$ we
say that $\sigma \in S_{n}$ is an $(j,n-j)${\bf -unshuffle} if
\[
\sigma(1) \leq\sigma(2) \leq \cdots \leq \sigma(j)
\hspace{.2in} \textrm{and} \hspace{.2in} \sigma(j+1) \leq
\sigma(j+2) \leq \cdots \leq \sigma(n).
\]
Readers familiar with shuffles will recognize unshuffles as their
inverses.

\begin{defn} \et \label{L-alg} An
\textbf{\textit{L}$_{\mathbf{\infty}}$-{\bf algebra}} is a graded
vector space $V = \bigoplus_{\ell = 0}^\infty V_{\ell}$
equipped with linear maps $l_{k} \maps V^{\otimes k} \rightarrow V$ 
for $k \ge 1$ with $\deg(l_{k}) = k-2$ such that:
\begin{eqnarray}
   l_{k}(x_{\sigma(1)}, \dots,x_{\sigma(k)}) =
   \chi(\sigma)l_{k}(x_{1}, \dots, x_{n})
\label{antisymmetry}
\end{eqnarray}
for all $\sigma \in S_{n}$ and $x_{1}, \dots, x_{n} \in V,$ and 
moreover, the following generalized form of the Jacobi identity
holds for $n \ge 0$:
\begin{eqnarray}
   \displaystyle{\sum_{i+j = n+1}
   \sum_{\sigma}
   \chi(\sigma)(-1)^{i(j-1)} l_{j}
   (l_{i}(x_{\sigma(1)}, \dots, x_{\sigma(i)}), x_{\sigma(i+1)},
   \ldots, x_{\sigma(n)}) =0,}
\label{megajacobi}
\end{eqnarray}
where the summation is taken over all $(i,n-i)$-unshuffles with $i
\geq 1.$
\end{defn}

In this definition the map $l_1$ makes $V$ into a chain complex,
since this map has degree $-1$ and Equation (\ref{megajacobi})
says its square is zero.  In what follows, we denote $l_1$ as $d$.
The map $l_{2}$ resembles a Lie bracket, since it is
skew-symmetric in the graded sense by Equation
(\ref{antisymmetry}). The map $l_3$ gives the Jacobiator and $l_4$
gives the Jacobiator identity.

To make this more precise, we make the following definition:

\begin{defn} \et A
\textbf{k-{\bf term} \textit{L}$_{\mathbf{\infty}}$-{\bf algebra}}
is an $L_{\infty}$-algebra $V$ with $V_{n} = 0$ for $n \geq k$.
\end{defn}

A $1$-term $L_{\infty}$-algebra is simply an ordinary Lie algebra,
where $l_{3} =0$ gives the Jacobi identity.  However, in a $2$-term
$L_{\infty}$-algebra, we no longer have $l_3 = 0$.
Instead, Equation (\ref{megajacobi}) says that the Jacobi identity
for $x,y,z \in V_0$ holds up to a term of the form
$dl_3(x,y,z)$.  We do, however, have $l_{4} = 0$, which provides
us with the coherence law that $l_{3}$ must satisfy.
It follows that a $2$-term $L_{\infty}$-algebra consists of:
\begin{itemize}
    \item vector spaces $V_{0}$ and
     $V_{1}$,

    \item a linear map $d\maps V_{1} \rightarrow V_{0},$

    \item bilinear maps $l_{2}\maps V_{i} \times V_{j}
     \rightarrow V_{i+j},$ where $0 \le i+j \le 1$,

    \item a trilinear map $l_{3}\maps V_{0} \times V_{0} \times
     V_{0} \rightarrow V_{1}$

\end{itemize}
satisfying a list of equations coming from Equations (\ref{antisymmetry})
and (\ref{megajacobi}) and the fact that $l_4 = 0$.
This list can be found in HDA6, but we will not need it here.

In fact, $2$-vector spaces are
equivalent to 2-term chain complexes of vector spaces: that
is, chain complexes of the form
\[                     V_1 \stackto{d} V_0 . \]
To obtain such a chain complex from a 2-vector space $L$, we let
$V_0$ be the space of objects of $L$.
However, $V_1$ is not the space of morphisms.  Instead, we define
the {\bf arrow part} $\vec{f}$ of a morphism $f \maps x \to y$ by
$$\vec{f} = f - i(s(f)), $$
and let $V_1$ be the space of these arrow parts.
The map $d \maps V_1 \to V_0$ is then just the target map $t
\maps \Mor(L) \to \Ob(L)$ restricted to $V_1 \subseteq \Mor(L)$.

To understand this construction a bit better, note
that given any morphism $f \maps x \to y$, its arrow
part is a morphism $\vec{f} \maps 0 \rightarrow y-x.$
Thus, taking the arrow part has the effect of `translating $f$
to the origin'.  We can always recover any morphism from its source together
with its arrow part, since $f = \vec{f} + i(s(f))$.   It follows that
any morphism $f \maps x \to y$ can be identified
with the ordered pair $(x, \vec{f})$ consisting of
its source and arrow part.  So, we have $\Mor(L) \iso V_0 \oplus V_1$.

We can actually recover the whole 2-vector space structure of $L$
from just the chain complex $d \maps V_1 \to V_0$.
To do this, we take:
\begin{eqnarray*}
  \Ob(L) & = & V_{0} \\
  \Mor(L) & = & V_{0} \oplus V_{1},
\end{eqnarray*}
with source, target and identity-assigning maps defined by:
\[
\begin{array}{ccl}
  s(x, \vec{f}) &=& x \\
  t(x, \vec{f}) &=& x + d\vec{f} \\
  i(x) &=& (x, 0)
\end{array}
\]
and with the composite of $f \maps x \to y$ and $g \maps y \to z$
defined by:
\[  g \circ f = (x, \vec{f} + \vec{g}).
\]
So, 2-vector spaces are equivalent to 2-term chain complexes.

Given this, it should not be surprising that Lie $2$-algebras
are equivalent to 2-term $L_\infty$-algebras.  Since we make
frequent use of this fact in the calculations to come, we recall
the details here.

Suppose $V$ is a 2-term $L_\infty$-algebra.   We obtain a 2-vector
space $L$ from the underlying chain complex of $V$ as above.
We continue by equipping $L$ with additional structure that
makes it a Lie $2$-algebra. It is sufficient to define the bracket
functor $[\cdot, \cdot] \maps L \times L \rightarrow L$ on a pair
of objects and on a pair of morphisms where one is an identity
morphism.  So, we set:
\[
\begin{array}{ccl}
[x,y] &=& l_2(x,y)       , \cr
[1_z,f] &=& (l_2(z,x), l_2(z,\vec{f})) , \cr
[f,1_z] &=& (l_2(x,z), l_2(\vec{f}, z))  ,
\end{array}
\]
where $f \maps x \to y$ is a morphism in $L$ and $z$ is an object.
Finally, we define the Jacobiator for $L$ in terms of its source and arrow
part as follows:
$$J_{x,y,z} = ([[x,y],z], l_{3}(x,y,z)).  $$
For a proof that $L$ defined this way is actually a Lie 2-algebra, see
HDA6.

In our calculations we shall often describe Lie 2-algebra homomorphisms
as homomorphisms between the corresponding $2$-term $L_{\infty}$-algebras:

\begin{defn} \et \label{Linftyhomo}
Let $V$ and $V'$ be $2$-term $L_{\infty}$-algebras. An
\textbf{\textit{L}$_{\mathbf{\infty}}$-{\bf homomorphism}}
$\phi \maps V \rightarrow V'$ consists of:
\begin{itemize}
\item a chain map $\phi \maps V \to V'$ consisting of linear maps
$\phi_0 \maps V_0 \to V'_0$  and
$\phi_1 \maps V_1 \to V'_1$,

\item a skew-symmetric bilinear map
$\phi_{2} \maps V_{0} \times V_{0} \to V_{1}'$,
\end{itemize}

such that the following equations hold for all $x,y,z \in V_0$
and $h \in V_{1}:$
\begin{eqnarray}
d (\phi_{2}(x,y)) = \phi_{0}(l_{2} (x,y)) - l_{2}(\phi_{0}(x), \phi_{0}(y))
\label{homo1}
\end{eqnarray}
\begin{eqnarray}
\phi_{2}(x,dh) = \phi_{1}(l_{2}(x,h)) - l_{2}(\phi_{0}(x), \phi_{1}(h))
\label{homo2}
\end{eqnarray}
\begin{equation}
\begin{array}{l}
\phi_1(l_3(x,y,z)) - l_3(\phi_0(x),\phi_0(y), \phi_0(z))
   = \\
   \phi_2(x,l_2(y,z))
   + \phi_2(y,l_2(z,x))
   + \phi_2(z,l_2(x,y)) \; +  \\
  {}   l_2(\phi_0(x),\phi_2(y,z))
  +    l_2(\phi_0(y),\phi_2(z,x))
  +    l_2(\phi_0(z),\phi_2(x,y))
\label{homo3}
\end{array}
\end{equation}
\end{defn}

\noindent
Equations (\ref{homo1}) and (\ref{homo2})
say that $\phi_{2}$ defines a chain homotopy from
$l_{2}(\phi(\cdot), \phi(\cdot))$ to $\phi(l_{2}(\cdot, \cdot))$,
where these are regarded as chain maps from $V \otimes V$ to $V'$.
Equation (\ref{homo3}) is just a chain complex
version of the commutative diagram in Definition
\ref{lie2algfunct}.   Furthermore, Definition 5.2 of Lada and Markl 
\cite{LM}, in the special case of 2-term $L_\infty$-algebras, reduces 
to the special case $l'_3=0$ of the definition above.

Let us sketch how to obtain
the Lie 2-algebra homomorphism $F$ corresponding to a given
$L_{\infty}$-homomorphism $\phi \maps V \to V'$. We define
the chain map $F \maps L \to L'$ in terms of $\phi$ using the
fact that objects of a 2-vector space are 0-chains in the corresponding
chain complex, while morphisms are pairs consisting of a 0-chain
and a 1-chain.  To make $F$ into a Lie 2-algebra homomorphism
we must equip it with a natural
transformation $F_{2}$ satisfying the conditions in Definition
\ref{lie2algfunct}.  In terms of its source and arrow parts, 
we define $F_2$ by
$$F_{2}(x,y) = (l_{2}(\phi_{0}(x), \phi_{0}(y)), \phi_{2}(x,y)).$$

We should also know how to compose $L_\infty$-homomorphisms.
We compose a pair of $L_{\infty}$-homomorphisms $\phi
\maps V \rightarrow V'$ and $\psi \maps V' \rightarrow V''$ by
letting the chain map $\psi \circ \phi \maps V \to V''$ be the
usual composite, while defining $(\psi \circ \phi)_{2}$ by:
\begin{equation}
 (\psi \circ \phi)_2(x,y) = \psi_2(\phi_0(x),\phi_0(y)) +
                                \psi_1(\phi_2(x,y)).
\label{composite.homo}
\end{equation}
This is just a chain complex version of how we compose
homomorphisms between Lie 2-algebras.  Note that the identity homomorphism
$1_V \maps V \to V$ has the identity chain map as its underlying
map, together with $(1_V)_2 = 0$.

We also have `2-homomorphisms' between homomorphisms:

\begin{defn} \et \label{Linfty2homo} Let $V$ and $V'$ be
2-term $L_\infty$-algebras and let $\phi, \psi \maps V \to V'$ be
$L_{\infty}$-homomorphisms.  An
\textbf{\textit{L}$_{\mathbf{\infty}}$-{\bf 2-homomorphism}}
$\tau \maps \phi \To \psi$ is a chain
homotopy $\tau$ from $\phi$ to $\psi$
such that the following equation holds for all $x,y \in
V_0$:
\begin{equation}
\phi_2(x,y) - \psi_2(x,y) =
l_{2}(\phi_{0}(x), \tau(y)) + l_{2}(\tau(x), \psi_{0}(y)) -
\tau ( l_{2}(x,y))
\label{2homo}
\end{equation}
\end{defn}

Given an $L_{\infty}$-2-homomorphism $\tau \maps \phi \To \psi$ between
$L_\infty$-homomorphisms $\phi, \psi \maps V \to V'$, there is a
corresponding Lie 2-algebra 2-homomorphism $\theta$ given by
\[
\theta(x) = (\phi_{0} (x), \tau(x))
\]

In HDA6, we showed:

\begin{prop} \et There is a strict $2$-category
{\bf 2TermL$_\mathbf\infty$} with $2$-term $L_{\infty}$-algebras
as objects, $L_\infty$-homomorphisms as morphisms, and
$L_\infty$-$2$-homomorphisms as $2$-morphisms.
\end{prop}

Using the equivalence between $2$-vector spaces and $2$-term chain
complexes, we established the equivalence between Lie $2$-algebras
and $2$-term $L_{\infty}$-algebras:

\begin{thm} \label{1-1'} \et The $2$-categories {\rm Lie$2$Alg}
and {\rm 2TermL$_{\infty}$} are $2$-equivalent.
\end{thm}

We use this result extensively in Section \ref{equivalence.section}.

\subsection{The Lie 2-Algebra $\g_k$}
\label{ghbar.section}

Thanks to the formula
\[
  d\vec{f} =  t(f) - s(f)   ,
\]
a 2-term $L_\infty$-algebra with vanishing differential corresponds
to a Lie 2-algebra for which the source of any morphism equals its
target.  In other words, the corresponding Lie $2$-algebra is `skeletal':

\begin{defn} \et
A category is {\bf skeletal} if isomorphic objects are always equal.
\end{defn}

Every category is equivalent to a skeletal one formed by choosing one
representative of each isomorphism class of objects \cite{Mac}.
As shown in HDA6, the same sort of thing is true for Lie $2$-algebras:

\begin{prop} \et \label{skeletal}
Every Lie $2$-algebra is equivalent, as an object of {\rm
Lie2Alg}, to a skeletal one.
\end{prop}

This result helps us classify Lie 2-algebras up to
equivalence.  We begin by reminding the reader of the relationship between
$L_{\infty}$-algebras and Lie algebra cohomology described in HDA6:

\begin{thm} \et \label{trivd}
There is a one-to-one correspondence between isomorphism classes
of $L_{\infty}$-algebras
consisting of only two nonzero terms $V_{0}$ and $V_{n}$ with
$d=0,$ and isomorphism classes of quadruples $(\g, V, \rho, [l_{n+2}])$
where $\g$ is a Lie algebra, $V$ is a vector space, $\rho$ is
a representation of $\g$ on $V$, and $[l_{n+2}]$ is an element
of the Lie algebra cohomology group $H^{n+2}(\g,V)$.
\end{thm}

\noindent
Here the representation $\rho$ comes from $\ell_2 \maps V_0 \times V_n
\to V_n$.

Because $L_{\infty}$-algebras are equivalent to Lie
2-algebras, which all have equivalent skeletal versions,
Theorem \ref{trivd} implies:

\begin{cor} \et \label{class} Up to equivalence, Lie 2-algebras
are classified by isomorphism classes of quadruples $(\g,\rho,V,[\ell_3])$
where:
\begin{itemize}
\item $\g$ is a Lie algebra,
\item $V$ is a vector space,
\item $\rho$ is a representation of $\g$ on $V$,
\item $[\ell_3]$ is an element of $H^3(\g,V)$.
\end{itemize}
\end{cor}

\noindent This classification lets us construct a
1-parameter family of Lie 2-algebras $\g_k$ for any simple real
Lie algebra $\g$:

\begin{example} \label{ghbar} \et Suppose $\g$ is a simple real
Lie algebra and $k \in \R$.  Then there is a skeletal Lie $2$-algebra
$\g_k$ given by taking $V_{0} = \g$, $V_{1} = \R$, $\rho$ the trivial
representation,  and $l_3(x,y,z) = k\langle x,[y,z]\rangle $.
\end{example}

\noindent
Here $\langle \cdot, \cdot \rangle$ is a suitably rescaled version
of the Killing form $\tr(\ad(\cdot)\ad(\cdot))$.
The precise rescaling factor will only
become important in Section \ref{loop.basic.section}.
The equation saying that $l_3$ is a $3$-cocycle is equivalent to the
equation saying that the left-invariant 3-form $\nu$ on $G$ with
$\nu(x,y,z) = \langle x,[y,z] \rangle$ is {\it closed}.

\subsection{The Lie 2-Algebra of a Fr\'echet Lie 2-Group}
\label{frechet.section}

Just as Lie groups have Lie algebras, `strict Lie 2-groups' have
`strict Lie 2-algebras'.  Strict Lie $2$-groups and Lie
$2$-algebras are categorified versions of Lie groups and Lie
algebras in which all laws hold `on the nose' as equations, rather
than up to isomorphism.    All the Lie 2-groups discussed in this
paper are strict.  However, most of them are infinite-dimensional 
`Fr\'echet' Lie 2-groups.

A {\bf Fr\'echet Lie group} is a Fr\'echet manifold \cite{Hamilton}
$G$ such that the multiplication map $m \maps G \times G \to G$ and
the inverse map $\inv \maps G \to G$ are smooth.  A {\bf homomorphism}
of Fr\'echet Lie groups is a group homomorphism that is also smooth.
For example, the space of smooth paths in a Lie group $G$ is a
Fr\'echet Lie group, and evaluation at a point defines a homomorphism
to $G$.  For more details we refer the reader to the survey article by
Milnor \cite{Milnor}, or Pressley and Segal's book on loop groups
\cite{PressleySegal}.

\begin{defn} \label{frechetlietwogroup}  \et
A {\bf strict Fr\'echet Lie $2$-group} $C$ is a
category such that
the set of objects $\Ob(C)$ and
the set of morphisms $\Mor(C)$
are both Fr\'echet Lie groups, and
the source and target maps $s, t \maps \Mor(C) \to
\Ob(L)$, the map $i \maps \Ob(C) \to \Mor(C)$ sending any object
to its identity morphism, and
the map $\circ \maps \Mor(C) \times_{\Ob(C)} \Mor(C) \to \Mor(C)$
sending any composable pair of morphisms to its composite
are all Fr\'echet Lie group homomorphisms.
\end{defn}

\noindent
Here $\Mor(C) \times_{\Ob(C)} \Mor(C)$ is the set of composable
pairs of morphisms, which we require to be a Fr\'echet Lie group.

Just as for ordinary Lie groups, taking the tangent space at the
identity of a Fr\'echet Lie group gives a Lie algebra.
Using this, it is not hard to see that strict Fr\'echet Lie 2-groups
give rise to Lie 2-algebras.  These Lie 2-algebras are actually
`strict':

\begin{defn} \et
A Lie 2-algebra is {\bf strict} if its Jacobiator is the identity.
\end{defn}

\noindent
This means that the map $l_3$ vanishes in the corresponding
$L_\infty$-algebra.  Alternatively:

\begin{prop} \et
A strict Lie 2-algebra is the same as a 2-vector space
$L$ such that $\Ob(L)$ and
$\Mor(L)$ are equipped with Lie algebra structures,
and the maps $s,t,i$ and $\circ$ are Lie algebra homomorphisms.
\end{prop}

\Proof -
A straightforward verification; see also HDA6.
\endofproof

\begin{prop} \label{lietwoalg.of.lietwogrp} \et
Given a strict Fr\'echet Lie $2$-group $\lietwogrp$, there is
a strict Lie $2$-algebra $\lietwoalg$ for which
$\Ob(\lietwoalg)$ and $\Mor(\lietwoalg)$ are the Lie algebras
of the Fr\'echet Lie groups $\Ob(\lietwogrp)$ and
$\Mor(\lietwogrp)$, respectively, and the maps
$s,t,i$ and $\circ$ are the differentials of the 
corresponding maps for $\lietwogrp$.
\end{prop}

\Proof  This is a generalization
of a result in HDA6 for ordinary Lie 2-groups,
which is straightforward to show directly.
\endofproof

In what follows all Fr\'echet Lie
2-groups are strict, so we omit the term `strict'.

\section{Review of Loop Groups}
\label{loop.section}

Next we give a brief review of loop groups and their central
extensions.  More details can be found in the canonical text
on the subject, written by Pressley and Segal \cite{PressleySegal}.

\subsection{Definitions and Basic Properties}
\label{loop.basic.section}

Let $G$ be a simply-connected compact simple Lie group.  We shall be
interested in the {\bf loop group} $\OG$ consisting of all smooth maps
from $[0,2\pi]$ to $G$ with $f(0) = f(2\pi) = 1$.  We make $\OG$ into
a group by pointwise multiplication of loops: $(fg)(\theta) =
f(\theta) g(\theta)$.  Equipped with its $C^\infty$ topology, $\OG$
naturally becomes an infinite-dimensional Fr\'echet manifold.  In fact
$\OG$ is a Fr\'echet Lie group, as defined in Section
\ref{frechet.section}.

As remarked by Pressley and Segal, the behavior of the group $\OG$
is ``untypical in its simplicity,'' since it turns out to behave
remarkably like a compact Lie group.  For example, it has an
exponential map that is locally one-to-one and onto, and it has a
well-understood highest weight theory of representations.  One
striking difference between $\OG$ and $G$, though, is the
existence of nontrivial central extensions of $\OG$ by the circle
$\TT$:
\begin{equation}
\label{eq: Kac--Moody extn}
1\to \TT \to \widehat{\Omega G}\stackrel{p}{\to} \Omega G\to 1 .
\end{equation}
It is important to understand that these extensions are
nontrivial, not merely in that they are classified by a nonzero
$2$-cocycle, but also \emph{topologically}.  In other words,
$\wOG$ is a nontrivial principal $\TT$-bundle over $\OG$ with the
property that $\wOG$ is a Fr\'echet Lie group, and $\TT$ sits
inside $\widehat{\Omega G}$ as a central subgroup in such a way
that the quotient $\widehat{\Omega G}/\TT$ can be identified with
$\Omega G$.  $\wOG$ is called the \textbf{Kac--Moody group}.

Associated to the central extension~\eqref{eq: Kac--Moody extn}
there is a central extension of Lie algebras:
\begin{equation}
\label{eq: Kac--Moody alg extn}
0 \to \u(1) \to \wOg \stackto{p} \Og \to 0
\end{equation}
Here $\Og$ is the Lie algebra of $\OG$,
consisting of all smooth maps $f \maps S^1\to \g$ such
that $f(0) = 0$.  The
bracket operation on $\Omega \mathfrak{g}$ is given by the
pointwise bracket of functions: thus $[f,g](\theta) =
[f(\theta),g(\theta)]$ if $f,g\in \Omega \mathfrak{g}$.
$\wOg$ is the simplest example of an affine Lie algebra.

The Lie algebra extension~\eqref{eq: Kac--Moody alg extn}
is simpler than the group extension~\eqref{eq: Kac--Moody extn}
in that it is determined up to isomorphism by a Lie algebra
$2$-cocycle $\omega(f,g)$, i.e.\ a skew bilinear map
$\omega\colon \Og\times \Og \to \R$ satisfying the
{\bf 2-cocycle condition}
\begin{equation}
\label{eq: 2-cocycle eqn}
\omega([f,g],h) + \omega([g,h],f) + \omega([h,f],g) = 0 .
\end{equation}
For $G$ as above we may assume the cocycle $\omega$ equal, up
to a scalar multiple, to the {\bf Kac--Moody $2$-cocycle}
\begin{equation}
\label{eq: Kac--Moody cocycle}
\omega(f,g) = 2\int^{2\pi}_0
\langle f(\theta),g'(\theta)\rangle\, d\theta
\end{equation}
where $\langle \cdot, \cdot \rangle$ is the basic inner product on
$\g$ divided by $4 \pi$.  Recall from Pressley and Segal
\cite{PressleySegal} that the basic inner product on $\g$ is the
unique invariant inner product $( \cdot, \cdot )$ with
$(h_{\theta},h_{\theta}) = 2$ where $h_{\theta}$ is the coroot
associated to the highest root $\theta$ of $\g$.
Thus, as a vector space $\wOg$ is isomorphic to $\Og\oplus \R$,
but the bracket is given by
$$
[(f,\a),(g,\b)] = ([f,g],\omega(f,g))
$$
Since $\omega$ is a skew form on $\Og$, it defines a left-invariant
2-form $\omega$ on $\OG$.  The cocycle
condition, Equation~\eqref{eq: 2-cocycle eqn},
says precisely that $\omega$ is closed.

For any  $k \in \R$, the cocycle $k\omega$ defines an extension of
Lie algebras
\[
0\to \u(1) \to \wOkg \to \Og \to 0
\]
where $\wOkg = \Og \oplus \u(1)$ with bracket defined in the same
way as for $\wOg$.  When $k$ is an integer, Pressley and Segal
\cite{PressleySegal} show that associated to this central
extension of Lie algebras is a central extension
\[
   1 \to \TT \to \wOkG \to \OG \to 1
\]
of Lie groups.  The integer $k$ is called the {\bf level} of the
central extension $\wOkG$.


\subsection{The Kac--Moody group $\wOkG$}
\label{Kac-Moody.section}

Several closely related explicit constructions of $\wOkG$ appear in the 
literature: first came the work Mickelsson \cite{Mickelsson}, then Murray 
\cite{Murray1}, then Brylinski--McLaughlin \cite{BM}, and more 
recently Murray--Stevenson \cite{MuSt}.  The last construction,
inspired by the work of Mickelsson, will be the 
most convenient for our purposes.  We shall use this 
to prove a result, Proposition \ref{conjugation}, that will be
crucial for constructing the 2-group $\PG$.

First, suppose that $\cG$ is any Fr\'echet Lie group.
Let $P_0\cG$ denote the space of smooth based paths in $\cG$:
\[   P_0\cG = \{ f \in C^\infty([0,2\pi], \cG) \colon \; f(0) = 1  \}
\]
$P_0\cG$ is a Fr\'echet Lie group under pointwise multiplication of
paths, whose Lie algebra is
\[   P_0L = \{ f \in C^\infty([0,2\pi], L) \colon \; f(0) = 0  \}
\]
where $L$ is the Lie algebra of $\cG$.  Furthermore, the map $\pi\colon
P_0\cG\to \cG$ which evaluates a path at its endpoint is a
homomorphism of Fr\'echet Lie groups.  The kernel
of $\pi$ is equal to
\[   \Omega\cG = \{ f \in C^\infty([0,2\pi], \cG) \colon \;
f(0) = f(2\pi) = 1  \, \}
\]
Thus, $\Omega \cG$ is a normal subgroup of $P_0 \cG$.  Note
that we are defining $\Omega\cG$ in a somewhat nonstandard way,
since its elements can be thought of as loops $f \maps S^1 \to \cG$
that are smooth everywhere except at the basepoint, where both left
and right derivatives exist to all orders, but need not agree.
However, we need this for the sequence
\[     1 \stackto{} \Omega\cG \stackto{} P_0 \cG \stackto{\pi} \cG
         \stackto{} 1 \]
to be exact.  Moreover, our $\Omega\cG$ is homotopy equivalent
to the usual group of smooth based loops in $\cG$.  
We give here a proof of this fact, due to Andrew Stacey \cite{Stacey}.  
For the purposes of this proof let $L_0\cG$ denote the group of
of smooth based maps from the circle $S^1$ into $\cG$, equipped with 
its $C^\infty$ topology.  $L_0\cG$ is a closed subgroup of $\Omega\cG$, 
so we have a continuous inclusion $i\colon L_0\cG \to \Omega\cG$.  
We construct a map $j\colon \Omega\cG \to 
L_0\cG$ as follows.  Let $\phi\colon [0,2\pi]\to [0,2\pi]$ 
be a smooth map of the interval $[0,2\pi]$ to itself, preserving the 
endpoints and with the property that all derivatives of $\phi$ 
vanish in some neighbourhood of the endpoints.  If $f\in \Omega\cG$ 
then $f\circ \phi \in L_0\cG$.  This defines a map 
$j\colon \Omega\cG \to L_0 \cG$, which 
is easily seen to be continuous.  To see that $j\circ i$ is homotopic to 
the identity map of $L_0\cG$, choose a linear homotopy 
$h_t$ from the identity map of the interval $[0,2\pi]$ to $\phi$.  
If $f\in L_0\cG$ then it is easy to check that $f\circ h_t$ 
is also in $L_0\cG$.  This defines a homotopy from $j\circ i$ 
to the identity.  A similar proof shows the composite $i\circ j$ 
is homotopic to the identity map of $\Omega\cG$.  

At present we are most interested in the case where $\cG = \OG$.
Then a point in $P_0\cG$ gives a map $f\maps [0,2\pi] \times
S^1 \to G$ with $f(0,\theta) = 1$ for all $\theta\in S^1$,
$f(t,0) = 1$ for all $t \in [0,2\pi]$.  Following 
\cite{MuSt}, we can see the map
$\kappa \colon P_0\OG\times P_0\OG \to \TT$ defined by
\begin{equation}
\label{eq: loop cocycle} \kappa(f,g) =
\exp\left(-2ik\int^{2\pi}_0\int^{2\pi}_0 \langle f(t)^{-1}f'(t),
g'(\theta) g(\theta)^{-1}\rangle\, d\theta \, dt\right)
\end{equation}
is a group $2$-cocycle.  This 2-cocycle $\kappa$ makes
$P_0\OG \times \TT$ into a group with the following product:
\[
  (f_1,z_1)\cdot (f_2,z_2) =
  \left(f_1f_2,z_1 z_2 \,\kappa(f_1,f_2)\right) .
\]
Let $N$ be the subset of $P_0\OG\times \TT$
consisting of pairs $(\c,z)$ such that $\c\colon [0,2\pi] \to \OG$
is a loop based at $1\in \OG$ and
\[
z = \exp\left(ik \int_{D_\c} \omega \right)
\]
where $D_\c$ is any disk in $\OG$ with $\gamma$ as its boundary.
It is easy to check that $N$
is a normal subgroup of the group $P_0\OG\times \TT$ with the product
defined as above.  To construct
$\wOkG$ we form the quotient group $(P_0\OG\times
\TT)/N$.  In \cite{MuSt} it is shown
that the resulting central extension is isomorphic to the
central extension of $\OG$ at level $k$.
So we have the commutative diagram
\begin{equation}
\label{centext}
\xymatrix{
P_0\OG \times \TT \ar[d] \ar[r] & \wOkG \ar[d] \\
P_0\OG \ar[r]^-{\pi} & \OG                                 }
\end{equation}
where the horizontal maps are quotient maps, the upper horizontal
map corresponding to the normal subgroup $N$, and the lower horizontal
map corresponding to the normal subgroup $\Omega^2G$ of $P_0\OG$.

Notice that the group of based paths $P_0G$ acts on $\OG$ by
conjugation.  The next proposition shows that this action
lifts to an action on $\wOkG$:

\begin{prop} \label{conjugation} \et
The action of $P_0G$ on $\OG$ by conjugation
lifts to a smooth action $\alpha$ of $P_0G$ on $\wOkG$, whose
differential gives an action $d\alpha$ of the Lie
algebra $P_0\g$ on the Lie algebra $\wOkg$ with
\[           d\alpha(p)(\ell,c) =
\big([p,\ell], \;
2k \int_0^{2\pi} \langle \ell(\theta), p'(\theta) \rangle \, d\theta \,) .
\]
for all $p \in P_0\g$ and all $(\ell,c) \in \Og \oplus \R
\iso \wOkg$.
\end{prop}

\Proof
To construct $\alpha$ it suffices to construct a smooth action of
$P_0G$ on $P_0\OG\times \TT$ that preserves the product
on this group and also preserves the normal subgroup $N$.
Let $p\colon [0,2\pi]\to G$ be an element of $P_0G$, so that
$p(0) = 1$.  Define the action of $p$ on a point $(f,z)\in
P_0\OG\times \TT$ to be
\[ p\cdot (f,z) =
\big( pfp^{-1}, \; z \exp(ik \int^{2\pi}_0 \beta_p(f(t)^{-1}f'(t))\,
dt)\; \big) \]
where $\beta_p$ is the left-invariant
$1$-form on $\OG$ corresponding to the linear map
$\beta_p\colon \Og\to \R$ given by:
$$
\beta_p(\xi) =
2 \int^{2\pi}_0 \langle \xi(\theta), p(\theta)^{-1} p'(\theta)\rangle \,
d\theta .
$$
for $\xi\in \Og$.
To check that this action preserves the product on
$P_0\OG\times \TT$, we have to show that
\[
\big(pf_1p^{-1}, \, z_1
\exp(ik\int^{2\pi}_0\beta_p(f_1(t)^{-1}f_1'(t))\, dt) \, \big)
\cdot
\big(pf_2p^{-1}, \, z_2
\exp(ik\int^{2\pi}_0\beta_p(f_2(t)^{-1}f_2'(t))\, dt) \, \big)  \]
\[
= \big(pf_1f_2p^{-1}, \,
z_1z_2 \kappa(f_1,f_2)
\exp(ik\int^{2\pi}_0 \beta_p((f_1f_2)(t)^{-1}(f_1f_2)'(t))\, dt) \,
\big).
\]
It therefore suffices to establish the identity
\begin{multline*}
\kappa(pf_1p^{-1},pf_2p^{-1}) =
\kappa(f_1,f_2)\exp \Big(ik\int^{2\pi}_0
\left(\beta_p((f_1f_2)(t)^{-1}(f_1f_2)'(t)) - \right. \\
\left.\beta_p(f_1(t)^{-1}f_1'(t)) - \beta_p(f_2(t)^{-1}f_2'(t))\right)
\, dt\Big) .
\end{multline*}
This is a straightforward computation that can safely
be left to the reader.

Next we check that the normal subgroup $N$ is preserved by the action
of $P_0G$.  For this we must show that if $(f,z)\in N$ then
$$
\big(pfp^{-1},\,z\exp(ik\int_0^{2\pi}\beta_p(f^{-1}f')\, dt)\big) \in N .
$$
Recall that $N$ consists
of pairs $(\c,z)$ such that $\c\in \Omega^2G$ and $z =
\exp(ik\int_{D_\c} \omega)$ where $D_\c$ is a disk in $\OG$ with
boundary $\c$.  Therefore we need to show that
$$
\exp\left(ik\int_{D_{p\c p^{-1}}}\omega \right) =
\exp\left(ik\int_{D_{\c}}\omega \right)
\exp\left(ik\int_0^{2\pi}\beta_p(\c^{-1}\c')dt \right) .
$$
This follows immediately from the identity
$$
\mathrm{Ad}(p)^*\omega = \omega + d\beta_p,
$$
which is easily established by direct computation.

Finally, we have to check the formula for $d\alpha$.
On passing to Lie algebras, diagram (\ref{centext})
gives rise to the following
commutative diagram of Lie algebras:
$$
\xymatrix{
P_0\Omega \mathfrak{g} \oplus \R \ar[r]^-{\overline{\mathrm{ev}}}
\ar[d] & \Omega \mathfrak{g} \oplus \R \ar[d] \\
P_0\Omega \mathfrak{g} \ar[r]^-{\mathrm{ev}} &
\Omega \mathfrak{g} }
$$
where $\overline{\mathrm{ev}}$ is the homomorphism
$(f,c)\mapsto (f(2\pi),c)$ for $f\in P_0\Omega \g$ and
$c\in \R$.  To calculate $d\alpha(p)(\ell,c)$ we compute
$\overline{\mathrm{ev}}(d\alpha(p)(\tilde{\ell},c))$ where
$\tilde{\ell}$ satisfies $\mathrm{ev}(\tilde{\ell}) = \ell$ (take,
for example, $\tilde{\ell}(t) = t\ell/2\pi$).  It is then
straightforward to compute that
$$
\overline{\mathrm{ev}}(d\alpha(p)(\tilde{\ell},c)) =
\big([p,\ell],2k\int^{2\pi}_0
\langle \ell(\theta),p'(\theta)\rangle\, d\theta\big).
$$
\hbox{\hskip 60ex} \endofproof

\section{The Lie 2-Group $\PG$}
\label{PG.section}

Having completed our review of Lie 2-algebras and loop groups,
we now study a Lie 2-group $\PG$ whose Lie
2-algebra $\Pg$ is equivalent to $\g_k$.   We begin in
Section \ref{PG.construction.section} by giving a construction of
$\PG$ in terms of the central extension $\wOkG$ of the loop group
of $G$.  This yields a description of $\Pg$ which we use later
to prove that this Lie 2-algebra is equivalent to $\g_k$.

Section \ref{topology.section} gives another viewpoint on $\PG$,
which goes a long way toward explaining the significance of this 2-group.
For this, we study the topological group $|\PG|$ formed by taking the
geometric realization of the nerve of $\PG$.

\subsection{Constructing $\PG$}
\label{PG.construction.section}

In Proposition \ref{conjugation} we saw that the action of
the path group $P_0G$ on the loop group $\OG$
by conjugation lifts to an action
$\alpha$ of $P_0G$ on the central extension $\wOkG$.
This allows us to define a Fr\'echet Lie group $P_0G
\ltimes \wOkG$ in which multiplication is given by:
\[
(p_1,\hat{\ell}_1)\cdot (p_2,\hat{\ell}_2) =
\big(\, p_1p_2, \, \hat{\ell_1} \alpha(p_1)(\hat{\ell}_2) \,\big) .
\]
This, in turn, allows us to construct the
2-group $\PG$ which plays the starring role in this paper:

\begin{prop} \label{PG.construction} \et
Suppose $G$ is a simply-connected compact simple Lie
group and $k \in \Z$.  Then there is a Fr\'echet Lie 2-group
$\PG$ for which:
\begin{itemize}
\item The Fr\'echet Lie group of objects $\Ob(\PG)$ is $P_0G$.
\item The Fr\'echet Lie group of morphisms $\Mor(\PG)$ is
$P_0G \ltimes \wOkG$.
\end{itemize}
\begin{itemize}
\item  The source and target maps $s,t \maps \Mor(\PG) \to \Ob(\PG)$
are given by:
\[
\begin{array}{ccl}
s(p,\hat{\ell}) &=& p  \\
t(p,\hat{\ell}) &=&  \partial(\hat{\ell}) p
\end{array}
\]
where $p \in P_0G$, $\hat{\ell} \in \wOkG$, and
$\partial \maps \wOkG \to P_0 G$ is the composite:
\[     \wOkG \to \OG \hookrightarrow P_0 G \, .\]
\item The identity-assigning map $i \maps \Ob(\PG) \to \Mor(\PG)$
is given by:
\[    i(p) = (p,1) . \]
\item  The composition map $\circ \maps \Mor(\PG) \times_{\Ob(\PG)}
\Mor(\PG) \to \Mor(\PG)$ is given by:
\[      (p_1,\hat{\ell}_1) \circ (p_2,\hat{\ell}_2) =
        (p_2, \hat{\ell}_1 \hat{\ell}_2)  \]
whenever $(p_1, \hat{\ell}_1), (p_2, \hat{\ell}_2)$ are composable
morphisms in $\PG$.
\end{itemize}
\end{prop}

\Proof  One can check directly that $s,t,i,\circ$ are Fr\'echet Lie
group homomorphisms and that these operations
make $\PG$ into a category.   Alternatively, one can check that
$(P_0 G, \wOkG, \alpha, \partial)$ is a crossed module
in the category of Fr\'echet manifolds.
This merely requires checking that
\begin{equation} \label{crossed.1}
   \partial( \alpha(p) (\hat{\ell})) =
p \, \partial(\hat{\ell})\, p^{-1}
\end{equation}
and
\begin{equation} \label{crossed.2}
    \alpha(\partial(\hat{\ell}_1)) (\hat{\ell_2}) =
\hat{\ell}_1 \hat{\ell}_2 {\hat{\ell}_1}^{-1} .
\end{equation}
Then one can use the fact that crossed modules in the
category of Fr\'echet manifolds are the same as Fr\'echet
Lie 2-groups (see for example HDA6).  \endofproof

We denote the Lie 2-algebra of $\PG$ by $\Pg$.  To prove
this Lie 2-algebra is equivalent to $\g_k$ in Section
\ref{equivalence.section}, we will use
an explicit description of its corresponding
$L_\infty$-algebra:

\begin{prop} \label{Pg} \et
The 2-term $L_\infty$-algebra $V$ corresponding to the
Lie 2-algebra $\Pg$ has:
\begin{itemize}
\item $V_{0} = P_0 \g$ and $V_{1} = \wOkg \iso \Og \oplus \R$,

\item $d\maps V_{1} \rightarrow V_{0}$ equal to
      the composite
     \[     \wOkg \to \Og \hookrightarrow P_0 \g \, ,\]

\item $l_2 \maps V_0 \times V_0 \to V_1$ given by
      the bracket in $P_0 \g$:
      \[       l_2(p_1,p_2) = [p_1,p_2],     \]
      and $l_2 \maps V_0 \times V_1 \to V_1$ given by the action
      $d\alpha$ of $P_0 \g$ on $\wOkg$, or explicitly:
      \[
      l_2(p, (\ell, c)) =
      \big([p,\ell],  \;
      2k\int_0^{2\pi} \langle \ell(\theta),
      p'(\theta) \rangle \, d\theta \; \big)
      \]
for all $p \in P_0\g$, $\ell \in \OG$ and $c \in \R$.

    \item $l_{3}\maps V_{0} \times V_{0} \times
     V_{0} \rightarrow V_{1}$ equal to zero.
\end{itemize}
\end{prop}

\Proof
This is a straightforward application of the correspondence
described in Section \ref{Linfty.section}.  The formula for
$l_2 \maps V_0 \times V_1 \to V_1$ comes from Proposition
\ref{conjugation}, while $\ell_3 = 0$
because the Lie 2-algebra $\Pg$ is strict.
\endofproof

\subsection{The Topology of $|\PG|$}
\label{topology.section}

In this section we construct an exact sequence of Fr\'echet Lie
2-groups:
\[
1\to \LG \stackto{\iota} \PG \stackto{\pi} G \to 1 \, ,
\]
where $G$ is considered as a Fr\'echet Lie $2$-group with only
identity morphisms.  Taking the geometric realization of the nerve,
we obtain this exact sequence of topological groups:
\[
1 \to |\LG| \stackto{|\iota|} |\PG| \stackto{|\pi|} G \to 1  \, .
\]
Note that $|G| = G$.  We then show that the topological group $|\LG|$
has the homotopy type of the Eilenberg--Mac Lane space $K(\Z,2)$.
Since $K(\Z,2)$ is also the classifying space $B\U(1)$, the above exact
sequence is a topological analogue of the exact sequence of Lie 2-algebras
describing how $\g_k$ is built from $\g$ and $\u(1)$:
\[   0 \to {\rm b}\u(1) \to \g_k \to \g \to 0 \, , \]
where ${\rm b}\u(1)$ is the Lie 2-algebra with a 0-dimensional
space of objects and $\u(1)$ as its space of morphisms.

The above exact sequence of topological groups exhibits $|\PG|$ as
the total space of a principal $K(\Z,2)$ bundle over $G$.  Bundles
of this sort are classified by their `Dixmier--Douady class',
which is an element of the integral third cohomology group of the base
space.  In the case at hand, this cohomology group is $H^3(G) \iso
\Z$, generated by the element we called $[\nu/2\pi]$ in the
Introduction.  We shall show that the Dixmier--Douady class of the
bundle $|\PG| \to G$ equals $k [\nu/2\pi]$.  Using this, we show that
for $k = \pm 1$, $|\PG|$ is a version of $\hat G$ --- the topological group
obtained from $G$ by killing its third homotopy group.

We start by defining a map $\pi \maps \PG \to G$ as follows.
We define $\pi$ on objects $p \in \PG$ as follows:
\[       \pi(p) = p(2 \pi) \in G . \]
In other words, $\pi$ applied to a based path in $G$ gives the
endpoint of this path.  We define $\pi$ on morphisms in the only
way possible, sending any morphism $(p,\hat{\ell}) \maps p \to
\partial(\hat{\ell}) p$
to the identity morphism on $\pi(p)$.  It is easy to see that
$\pi$ is a {\bf strict homomorphism} of Fr\'echet Lie 2-groups:
in other words, a map that strictly preserves all the Fr\'echet
Lie 2-group structure.  Moreover, it is easy to see that
$\pi$ is onto both for objects and morphisms.

Next, we define the Fr\'echet Lie 2-group $\LG$ to be the
{\bf strict kernel} of $\pi$.   In other words, the objects of $\LG$ are
objects of $\PG$ that are mapped to $1$ by $\pi$,
and similarly for the morphisms of $\LG$, while the source,
target, identity-assigning and composition maps for $\LG$ are
just restrictions of those for $\PG$.  So:
\begin{itemize}
\item the Fr\'echet Lie group of objects $\Ob(\LG)$ is $\OG$,
\item the Fr\'echet Lie group of morphisms $\Mor(\LG)$ is
$\OG \ltimes \wOkG$,
\end{itemize}
where the semidirect product is formed using the action
$\alpha$ restricted to $\OG$.  Moreover, the formulas for
$s,t,i,\circ$ are just as in Proposition \ref{PG.construction},
but with loops replacing paths.

It is easy to see that the inclusion $\iota \maps \LG \to \PG$
is a strict homomorphism of Fr\'echet Lie 2-groups.
We thus obtain:

\begin{prop} \label{exact1} \et
The sequence of strict Fr\'echet 2-group homomorphisms
\[
1\to \LG \stackto{\iota} \PG \stackto{\pi} G \to 1
\]
is {\bf strictly exact}, meaning that the image of each arrow
is equal to the kernel of the next, both on objects and on morphisms.
\end{prop}

Any Fr\'echet Lie 2-group $C$ is, among other things, a {\bf
topological category}: a category where the sets $\Ob(C)$ and
$\Mor(C)$ are topological spaces and the source, target,
identity-assigning and composition maps are continuous.  
There is a standard procedure for taking the `nerve' of a
topological category and obtaining a simplicial space.  One
can then take the `geometric realization' of any simplicial
space, obtaining a topological space.  We use $|C|$ to denote the
geometric realization of the nerve of a topological category $C$.
If $C$ is in fact a topological 2-group --- for example a
Fr\'echet Lie 2-group --- then $|C|$ naturally becomes a
topological group \cite{Segal}.

Applying the functor $|\cdot|$ to the exact sequence in
Proposition \ref{exact1}, we obtain this result, which
implies Theorem 2:

\begin{thm} \label{exact2} \et
The sequence of topological groups
\[
1\to |\LG| \stackto{|\iota|} |\PG| \stackto{|\pi|} G \to 1
\]
is exact, and $|\LG|$ has the homotopy type of $K(\Z,2)$.
Thus, $|\PG|$ is the total space of a $K(\Z,2)$ bundle over $G$.
The Dixmier--Douady class of this bundle is $k [\nu/2\pi] \in
H^3(G)$.   Moreover, $\PG$ is $\hat G$ when $k = \pm 1$.
\end{thm}

\Proof
It is easy to see directly that the functor $|\cdot|$ carries strictly
exact sequences of topological 2-groups to exact sequences of
topological groups.  To show that $|\LG|$ is a $K(\Z,2)$, we prove
there is a strictly exact sequence of Fr\'echet Lie 2-groups
\begin{equation}
\label{exact.sequence}
1 \to \U(1) \to \T\wOkG \to \LG \to 1 \, .
\end{equation}
Here $\U(1)$ is regarded as a Fr\'echet Lie 2-group with only identity
morphisms, while $\T\wOkG$ is the Fr\'echet Lie 2-group with $\wOkG$
as its Fr\'echet Lie group of objects and precisely one morphism from
any object to any other.  In general:

\begin{lem} \label{T} \et
For any Fr\'echet Lie group
$\cG$, there is a Fr\'echet Lie 2-group $\T\cG$ with:
\begin{itemize}
\item $\cG$ as its Fr\'echet Lie group of objects,
\item $\cG \ltimes \cG$ as its Fr\'echet Lie group of morphisms,
where the semidirect product is defined using the conjugation
action of $\cG$ on itself,
\end{itemize}
and:
\begin{itemize}
\item source and target maps given by $s(g,g') = g$, $t(g,g') = gg'$,
\item identity-assigning map given by $i(g) = (g,1)$,
\item composition map given by $(g_1,g_1') \circ (g_2,g_2') =
(g_2,g_1'g_2')$ whenever $(g_1,g_1')$, $(g_2,g_2')$ are composable
morphisms in $\T\cG$.
\end{itemize}
\end{lem}

\Proof
It is straightforward to check that
this gives a Fr\'echet Lie 2-group.  Note that
$\T\cG$ has $\cG$ as objects and
one morphism from any object to any other.
\hbox{\hskip 60ex} \endofproof

In fact, Segal \cite{Segal} has already introduced $\T\cG$ under the
name $\overline{\cG}$, treating it as a topological category.  He
proved that $|\T\cG|$ is contractible.  In fact, he exhibited
$|\T\cG|$ as a model of $E\cG$, the total space of the universal bundle
over the classifying space $B\cG$ of $\cG$.
Therefore, applying the functor $| \cdot |$ to the exact
sequence (\ref{exact.sequence}), we obtain this short exact sequence of
topological groups:
\[
1 \to \U(1) \to E\wOkG \to |\LG| \to 1 \, .
\]
Since $E\wOkG$ is contractible, it follows
that $|\LG| \iso E\wOG/\U(1)$ has the homotopy
type of $B\U(1) \simeq K(\Z,2)$.

To see that $|\pi|\colon |\mathcal{P}_kG|\to G$ is
locally trivial, let $U\subset G$ be open; then the inverse image
$|\mathcal{P}_kG|_U = |\pi|^{-1}(U)$ is the geometric realization
of the nerve of the topological category with objects $P_0G|_U$
and morphisms $P_0G|_U\ltimes \hat{\Omega}_k G$, where 
$P_0G|_U$ denotes the inverse image of $U\subset G$ 
under the projection $P_0G\to G$.  Since $P_0G\to
G$ is a locally trivial principal bundle with structure group
$\Omega G$, we can find a homeomorphism
$$
P_0G|_U \stackrel{\sim}{\rightarrow} U\times \Omega G.
$$
This homeomorphism induces an isomorphism of categories
$$
\xymatrix{ P_0G|_U \ltimes \hat{\Omega}_k G \ar@<0.5ex>[d]
\ar@<-0.5ex>[d] \ar[r] & U\times \Omega G\ltimes \hat{\Omega}_k G
\ar@<0.5ex>[d]
\ar@<-0.5ex>[d] \\
P_0G|_U \ar[r] & U\times \Omega G }
$$
Since the functor $|\cdot |$ preserves products, the 
geometric realization of the topological category with
objects $U\times \Omega G$ and morphisms $U\times \Omega G\ltimes
\hat{\Omega}_k G$ is just $U\times |\mathcal{L}_k G|$.  Therefore the 
isomorphism of categories above induces, on taking geometric 
realisations, a
homeomorphism
$$
|\mathcal{P}_k G|_U \stackrel{\sim}{\rightarrow} U\times
|\mathcal{L}_k G|
$$
commuting with the projections to $U\subset G$.  
Moreover this homeomorphism is clearly $|\mathcal{L}_k G|$
equivariant, thus $|\mathcal{P}_kG|\to G$ is a locally trivial
principal bundle.
Like any such bundle, this is the pullback of the universal
principal $K(\Z,2)$ bundle $p \maps EK(\Z,2) \to BK(\Z,2)$ along
some map $f \maps G \to BK(\Z,2)$, giving a commutative diagram of
spaces:
\[
\begin{diagram}
\node{|\LG|} \arrow{e,t}{|\iota|} \arrow{s,l}{\sim}
\node{|\PG|} \arrow{e,t}{|\pi|} \arrow{s,r}{p^\ast f}
\node{G} \arrow{s,r}{f}                   \\
\node{K(\Z,2)} \arrow{e,t}{i}
\node{EK(\Z,2)} \arrow{e,t}{p}
\node{BK(\Z,2)}
\end{diagram}
\]
Indeed, such bundles are classified up to isomorphism by the homotopy class
of $f$.  Since $BK(\Z,2) \simeq K(\Z,3)$, this homotopy class is determined
by the Dixmier--Douady class $f^\ast \kappa$, where $\kappa$ is the
generator of $H^3(K(\Z,3)) \iso \Z$.  The next order of business is
to show that $f^\ast \kappa = k[\nu/2\pi]$.

For this, it suffices to show that $f$ maps the generator of
$\pi_3(G) \iso \Z$ to $k$ times the generator of $\pi_3(K(\Z,3))
\iso \Z$.   Consider this bit of the long exact sequences of homotopy
groups coming from the above diagram:
\[
\begin{diagram}
\node{\pi_3(G)} \arrow{e,t}{\partial} \arrow{s,l}{\pi_3(f)}
\node{\pi_2(|\LG|)} \arrow{s,r}{\iso}   \\
\node{\pi_3(K(\Z,3))} \arrow{e,t}{\partial'} \node{\pi_2(K(\Z,2))}
\end{diagram}
\]
Since the connecting homomorphism $\partial'$ and the
map from $\pi_2(|\LG|)$ to $\pi_2(K(\Z,2))$ are isomorphisms,
we can treat these as the identity by a suitable choice of
generators.  Thus, to show that $\pi_3(f)$ is multiplication by $k$
it suffices to show this for the connecting homomorphism $\partial$.

To do so, consider this commuting diagram of Frech\'et Lie 2-groups:
\[
\begin{diagram}
\node{\OG} \arrow{e,t}{\iota} \arrow{s,l}{i}
\node{P_0 G} \arrow{e,t}{\pi} \arrow{s,l}{i'}
\node{G} \arrow{s,r}{1}
              \\
\node{\LG} \arrow{e,t}{\iota}
\node{\PG}  \arrow{e,t}{\pi}
\node{G}
\end{diagram}
\]
Here we regard the groups on top as 2-groups with only identity
morphisms; the downwards-pointing arrows include these in
the 2-groups on the bottom row.  Applying the functor $|\cdot|$,
we obtain a diagram where each row is a principal bundle:
\[
\begin{diagram}
\node{\OG} \arrow{e,t}{|\iota|} \arrow{s,l}{|i|}
\node{P_0 G} \arrow{e,t}{|\pi|} \arrow{s,l}{|i'|}
\node{G} \arrow{s,r}{1}
              \\
\node{|\LG|} \arrow{e,t}{|\iota|}
\node{|\PG|}  \arrow{e,t}{|\pi|}
\node{G}
\end{diagram}
\]
Taking long exact sequences of homotopy groups, this gives:
\[
\begin{diagram}
\node{\pi_3(G)} \arrow{e,t}{1} \arrow{s,l}{1}
\node{\pi_2(\OG)} \arrow{s,r}{\pi_2(|i|)}
              \\
\node{\pi_3(G)} \arrow{e,t}{\partial}
\node{\pi_2(|\LG|)}
\end{diagram}
\]
Thus, to show that $\partial$ is multiplication by $k$
it suffices to show this for $\pi_2(|i|)$.

For this, we consider yet another commuting diagram of
Frech\'et Lie 2-groups:
\[
\begin{diagram}
\node{\U(1)} \arrow{e} \arrow{s}
\node{\wOkG} \arrow{e} \arrow{s}
\node{\OG}             \arrow{s,l}{i}
                          \\
\node{\U(1)}   \arrow{e}
\node{\T\wOkG} \arrow{e}
\node{\LG}
\end{diagram}
\]
Applying $|\cdot|$, we obtain a diagram where each row is a
principal $\U(1)$ bundle:
\[
\begin{diagram}
\node{\U(1)} \arrow{e} \arrow{s}
\node{\wOkG} \arrow{e} \arrow{s}
\node{\OG}             \arrow{s,r}{|i|}
                          \\
\node{\U(1)}   \arrow{e}
\node{|\T\wOkG|} \arrow{e}
\node{|\LG| \simeq K(\Z,2)}
\end{diagram}
\]
Recall that the bottom row is the universal principal $\U(1)$ bundle.
The arrow $|i|$ is the classifying map for the $\U(1)$ bundle $\wOkG \to \OG$.
The Chern class of this bundle
is $k$ times the generator of $H^2(\OG)$ (see for 
instance \cite{PressleySegal}), so $\pi_2(|i|)$ must map the
generator of $\pi_2(\OG)$ to $k$ times the generator of $\pi_2(K(\Z,2))$.

Finally, let us show that $|\PG|$ is $\hat{G}$ when $k = \pm 1$.
For this, it suffices to show that when $k = \pm 1$, the map
$|\pi| \maps |\PG| \to G$ induces isomorphisms on all homotopy
groups except the third, and that $\pi_3(|\PG|) = 0$.  For this
we examine the long exact sequence:
\[        \cdots \stackto{}
    \pi_n(|\LG|) \stackto{} \pi_n(|\PG|)
                 \stackto{} \pi_n(G)
                 \stackto{\partial} \pi_{n-1}(|\LG|)
          \stackto{} \cdots  \, .
\]
Since $|\LG| \simeq K(\Z,2)$, its homotopy groups vanish
except for $\pi_2(|\LG|) \iso \Z$, so $|\pi|$ induces an
isomorphism on $\pi_n$ except possibly for $n = 2,3$.
In this portion of the long exact sequence we have
\[
                    0
                \stackto{} \pi_3(|\PG|)
                 \stackto{} \Z
                 \stackto{k} \Z
                 \stackto{} \pi_2(|\PG|)
                 \stackto{} 0
\]
so $\pi_3(|\PG|) \iso 0$ unless $k = 0$, and
$\pi_2(|\PG|) \iso \Z/k\Z$, so $\pi_2(|\PG|) \iso \pi_2(G) \iso 0$
when $k = \pm 1$.
\endofproof

\section{The Equivalence Between $\Pg$ and $\g_k$}
\label{equivalence.section}

In this section we prove our main result, which implies
Theorem 1:

\begin{thm} \et \label{bigthm}
There is a strictly exact sequence of Lie 2-algebra homomorphisms
\[
    0 \to \T\Og \stackto{\lambda}
     \Pg \stackto{\phi}  \g_k  \to 0  \]
where $\T\Og$ is equivalent to the trivial Lie 2-algebra
and $\phi$ is an equivalence of Lie 2-algebras.
\end{thm}

\noindent
Recall that by `strictly exact' we mean that both on the vector
spaces of objects and the vector spaces of morphisms,
the image of each map is the kernel of the next.

We prove this result in a series of lemmas.  We begin by describing
$\T{\Og}$ and showing that it is equivalent to the trivial Lie
2-algebra.  Recall that in Lemma \ref{T} we constructed for any
Fr\'echet Lie group $\cG$ a Fr\'echet Lie 2-group $\T\cG$ with $\cG$
as its group of objects and precisely one morphism from any object
to any other.  We saw that the space $|\T\cG|$ is contractible; this is
a topological reflection of the fact that $\T\cG$ is equivalent to the
trivial Lie 2-group.  Now we need the Lie algebra analogue of this
construction:

\begin{lem} \label{TL} \et Given a Lie algebra $L$, there is
a 2-term $L_{\infty}$-algebra $V$ for which:

\begin{itemize}
    \item $V_{0} = L$ and $V_{1} = L$,

    \item $d\maps V_{1} \rightarrow V_{0}$ is the identity,

    \item $l_2 \maps V_0 \times V_0 \to V_1$ and
    $l_2 \maps V_0 \times V_1 \to V_1$
    are given by the bracket in $L$,

    \item $l_{3}\maps V_{0} \times V_{0} \times
     V_{0} \rightarrow V_{1}$ is equal to zero.
\end{itemize}
We call the corresponding strict Lie 2-algebra $\T{L}$.
\end{lem}

\Proof
Straightforward.
\endofproof

\begin{lem} \label{trivial} \et
For any Lie algebra $L$, the Lie 2-algebra $\T{L}$ is equivalent to the
trivial Lie 2-algebra.  That is, $\T{L} \simeq 0$.
\end{lem}

\Proof There is a unique homomorphism $\beta \maps \T{L} \to 0$
and a unique homomorphism $\gamma \maps 0 \to \T{L}$.  Clearly
$\beta \circ \gamma$ equals the identity.  The composite $\gamma \circ
\beta$ has:
\[
\begin{array}{lccl}
    (\gamma \circ \beta)_0 \maps &x &\mapsto& 0
    \\
    (\gamma \circ \beta)_1 \maps &x &\mapsto& 0
    \\
    (\gamma \circ \beta)_2 \maps &(x_1,x_2) &\mapsto& 0 \, ,
\end{array}
\]
while the identity homomorphism from $\T{L}$ to itself has:
\[
\begin{array}{lccl}
    \mathrm{id}_0 \maps &x &\mapsto& x
    \\
    \mathrm{id}_1 \maps &x &\mapsto& x
    \\
    \mathrm{id}_2 \maps &(x_1,x_2) &\mapsto& 0
    \,.
\end{array}
\]
There is a 2-isomorphism
\[
    \tau \maps \gamma \circ \beta \stackTo{\sim} \mathrm {id}
\]
given by
\[
    \tau\of{x} = x
    \,,
\]
where the $x$ on the left is in $V_0$ and that on the right in
$V_1$, but of course $V_0 = V_1$ here. \endofproof

We continue by defining the Lie 2-algebra homomorphism $\Pg
\stackto{\phi} \g_k$.

\begin{lem} \label{phi} \et
There exists a Lie 2-algebra homomorphism
\[
    \phi \maps \Pg \to \g_k
\]
which we describe in terms of its corresponding $L_{\infty}$-homomorphism:
\[
\begin{array}{ccl}
\phi_0 (p) &=& p(2\pi) \\
           &&          \\
\phi_1 (\ell, c) &=& c \\
           &&          \\
\phi_2 (p_1, p_2) &=&
       k \displaystyle{\int_0^{2\pi}} \left( \bracket{p_2}{p_1'}
      - \bracket{p_2'}{p_1} \right) \, d\theta
\end{array}
\]
where $p, p_1, p_2 \in P_0 \g$ and
$(\ell, c) \in \Og \oplus \R \iso \wOkg$.
\end{lem}

\noindent Before beginning, note that the quantity
\[
    \displaystyle{\int_0^{2\pi} \! \left( \bracket{p_2}{p_1'} -
     \bracket{p_2'}{p_1} \right) \, d\theta
     =
    2  \int_0^{2\pi} \! \bracket{p_2}{p_1'} \; d\theta}
    \; - \;
     \bracket{p_2(2\pi)}{p_1 (2\pi)}
\]
is skew-symmetric, but not in general equal to
\[
    \displaystyle{2\int_0^{2\pi}
       \bracket{p_2}{p_1'} \; d\theta}
\]
due to the boundary term.  However, these quantities are equal
when either $p_1$ or $p_2$ is a loop.

\vskip 1em

\Proof We must check that $\phi$ satisfies the conditions of Definition
\ref{Linftyhomo}.  First we show that $\phi$ is a chain map.
That is, we show that $\phi_0$ and $\phi_1$ preserve the differentials:
$$\xymatrix{
    \wOkg \ar[rr]^{d}
    \ar[dd]_{\phi_{1}} &&
    P_0 \g  \ar[dd]^{\phi_{0}} \\ \\
    \R \ar[rr]^{d'} &&
    \g
}$$
where $d$ is the composite given in Proposition \ref{Pg}, and $d'=0$
since $\g_k$ is skeletal.  This square commutes since $\phi_0$ is
also zero.

We continue by verifying conditions (\ref{homo1}) - (\ref{homo3}) of
Definition \ref{Linftyhomo}.  The bracket on objects is preserved
on the nose, which implies that the right-hand side of (\ref{homo1}) is
zero.  This is consistent with the fact that the differential in
the $L_\infty$-algebra for $\g_k$ is zero, which implies that the
left-hand side of (\ref{homo1}) is also zero.

The right-hand side of (\ref{homo2}) is given by:
\begin{eqnarray*}
      \phi_1( l_2(p,(\ell,c))
      -
      l_2 (
        \phi_0 (p),
        \phi_1 (\ell,c)
      ) &=&
      \phi_1 \Big( [p,\ell],
      2k\int \bracket{\ell}{p'} \; d\theta \Big) \
      -
      \underbrace{l_2(p(2\pi),c)}_{=0} \cr
      &=&
      2k \int \bracket{\ell}{p'}
      \; d\theta .
    \end{eqnarray*}
This matches the left-hand side of (\ref{homo2}), namely:
\begin{eqnarray*}
    \phi_2 (p, d (\ell,c)  )
    &=&
    \phi_2 (p, \ell) \\
    &=&
    k
    \int
      (\bracket{\ell}{p'}
      -
      \bracket{\ell'}{p} ) \; d\theta
\nonumber \cr
    &=&
    2k \int
      \bracket{\ell}{p'}
    \; d \theta
\end{eqnarray*}
Note that no boundary term appears here
since one of the arguments is a loop.

Finally, we check condition (\ref{homo3}).  Four terms in
this equation vanish because $l_3 = 0$ in $\Pg$ and
$l_2 = 0$ in $\g_k$.  We are left needing to show
\[
-l_3 (\phi_0 (p_1),\phi_0 (p_2), \phi_0 (p_3)) =
\phi_2 (p_1,l_2 (p_2,p_3)) + \phi_2 (p_2,l_2 (p_3,p_1)) +
\phi_2 (p_3,l_2 (p_1,p_2)) \, .
\]
The left-hand side here equals
$-k \langle p_1(2\pi), [p_2(2\pi), p_3(2\pi)]\rangle $.
The right-hand side equals:
\[
\begin{array}{ccl}
&& \phi_2(p_1, l_2(p_2,p_3)) \;\; + \;\; {\rm cyclic\; permutations}   \\
\\
&=&
k \displaystyle{\int \left( \bracket{[p_2,p_3]}{p_1'}
-\bracket{[p_2,p_3]'}{p_1} \right) \, d\theta \;\; + \; {\rm cyclic\; perms.} }
\\
\\
&=&
k \displaystyle{ \int \left( \bracket{[p_2,p_3]}{p'_1}
             +\bracket{[p'_2,p_3]}{p_1}
             -\bracket{[p_2,p'_3]}{p_1} \right) \, d\theta \;\; + \;\;
{\rm cyclic \; perms.} } \\
\end{array}
\]
Using the antisymmetry of $\langle \cdot, [\cdot, \cdot] \rangle$, this
becomes:
\[ k \int \left( \bracket{p'_1}{[p_2,p_3]}
                -\bracket{p'_2}{[p_3,p_1]}
                -\bracket{p'_3}{[p_1,p_2]} \right)
                 \, d\theta \;\; + \;\; {\rm cyclic \; perms.} \\
\]
The first two terms cancel when we add all their cyclic permutations,
so we are left with all three cyclic permutations of the last term:
\[ -k \int \left( \bracket{p'_1}{[p_2,p_3]}
                +\bracket{p'_2}{[p_3,p_1]}
                +\bracket{p'_3}{[p_1,p_2]} \right) \, d\theta \, . \]
If we apply integration by parts to the first term, we get:
\[ -k \int \left(-\bracket{p_1}{[p'_2,p_3]} - \bracket{p_1}{[p_2,p'_3]}
                +\bracket{p'_2}{[p_3,p_1]}
                +\bracket{p'_3}{[p_1,p_2]} \right) \; d\theta \; - \;  \]
\[ k \bracket{p_1(2\pi)}{[p_2(2\pi),p_3(2\pi)]} \, .\]
By the antisymmetry of $\langle \cdot, [\cdot, \cdot] \rangle$,
the four terms in the integral cancel, leaving just
$-k \langle p_1(2\pi), [p_2(2\pi), p_3(2\pi)] \rangle$, as desired.
\endofproof

Next we show that the strict kernel of
$\phi \maps \Pg \to \g_k$ is $\T\Omega\g$:

\begin{lem} \et There is a Lie 2-algebra homomorphism
\[   \lambda \maps \T{\Omega\g} \to \Pg, \]
that is one-to-one both on objects and on morphisms, and
whose range is precisely the kernel of
$\phi \maps \Pg \to \g_k$, both on objects and on morphisms.
\end{lem}

\Proof
Glancing at the formula for $\phi$ in Lemma \ref{phi}, we
see that the kernel of $\phi_0$ and the kernel of
$\phi_1$ are both $\Og$.  We see from Lemma \ref{TL} that
these are precisely the spaces $V_0$ and $V_1$ in the
2-term $L_\infty$-algebra $V$ corresponding to $\T{\Og}$.
The differential $d \maps \ker(\phi_1) \to \ker(\phi_0)$
inherited from $\T\Og$ also matches that in $V$: it is
the identity map on $\Og$.

Thus, we obtain an inclusion of 2-vector spaces
$\lambda \maps  \T{\Omega\g} \to \Pg$.
This uniquely extends to a Lie $2$-algebra homomorphism, which we
describe in terms of its corresponding $L_{\infty}$-homomorphism:
{
\begin{eqnarray*}
\lambda_0 (\ell) &=& \ell \cr
&&                         \cr
\lambda_1 (\ell) &=& (\ell, 0) \cr
\lambda_2 (\ell_1, \ell_2) &=&
\displaystyle{\big(0, 2k
             \int_0^{2\pi} \bracket{\ell_1}{\ell_2'} \, d\theta \big)}
\end{eqnarray*}
}
where $\ell, \ell_1, \ell_2 \in \Omega\g$, and the zero in the last
line denotes the zero loop.

To prove this, we must show that the conditions of Definition \ref{Linftyhomo}
are satisfied.  We first check that $\lambda$ is a chain map, i.e., this
square commutes:
$$\xymatrix{
    \Omega\g \ar[rr]^{d}
    \ar[dd]_{\lambda_{1}} &&
    \Omega\g   \ar[dd]^{\lambda_{0}} \\ \\
     \wOkg \ar[rr]^{d'} &&
    P_0 \g
}$$
where $d$ is the identity and $d'$ is the composite given in
Proposition \ref{Pg}.  To see this, note that
$d'(\lambda_1 (\ell)) = d'(\ell,0) = \ell$ and
$\lambda_0 (d(\ell)) = \lambda_0 (\ell) = \ell$.

We continue by verifying conditions (\ref{homo1}) - (\ref{homo3}) of
Definition \ref{Linftyhomo}.
The bracket on the space $V_0$ is strictly preserved
by $\lambda_0$, which implies that the right-hand side of
(\ref{homo1}) is zero.  It remains to show that the left-hand side,
$d'(\lambda_2 (\ell_1, \ell_2))$, is also zero.
Indeed, we have:
\[
d'(\lambda_2 (\ell_1, \ell_2)) =
d'\left(0, 2k \int \bracket{\ell_1}{\ell_2'} \; d\theta \right) =
0 \, .
\]

Next we check property (\ref{homo2}).  On the
right-hand side, we have:
\begin{eqnarray*}
        \lambda_1 ( l_2 (\ell_1,\ell_2) )
        -
        l_2 (\lambda_0 (\ell_1) , \lambda_1 (\ell_2))
        &=&
        ([\ell_1, \ell_2],0)
        -
        \big( [\ell_1, \ell_2],
        -2k \int \bracket{\ell_1}{\ell_2 '} \, d\theta \big) \\
        &=&
        \big(0, 2k  \int \bracket{\ell_1}{\ell_2 '} \, d\theta \big) \, .
\end{eqnarray*}
On the left-hand side, we have:
$$\lambda_2(\ell_1, d(\ell_2)) = \lambda_2 (\ell_1, \ell_2)
= \big(0, 2k \int \bracket{\ell_1}{\ell_2 '} \; d\theta \big)$$
Note that this also shows that given the chain map defined by
$\lambda_0$ and $\lambda_1$, the function $\lambda_2$ that extends
this chain map to an $L_\infty$-homomorphisms is unquely
fixed by condition (\ref{homo2}).

Finally, we show that $\lambda_2$ satisfies condition (\ref{homo3}).
The two terms involving $l_3$ vanish since $\lambda$ is a map between two
strict Lie $2$-algebras.  The three terms of the form
$l_2 ( \lambda_0 (\cdot), \lambda_2 (\cdot, \cdot))$ vanish because
the image of $\lambda_2$ lies in the center of $\wOkg$.
It thus remains to show that
$$
\lambda_2(\ell_1, l_2(\ell_2, \ell_3)) +
\lambda_2 (\ell_2, l_2 (\ell_3, \ell_1))
+ \lambda_2(\ell_3, l_2(\ell_1, \ell_2)) = 0. $$
This is just the cocycle property of the Kac--Moody cocycle,
Equation (\ref{eq: 2-cocycle eqn}).
\hbox{\hskip 60ex} \endofproof

Next we check the exactness of the sequence
\[    0 \to \T{\Omega\g} \stackto{\lambda}
            \Pg \stackto{\phi} \g_k  \to 0  \]
at the middle point.  Before doing so, we recall the formulas for
the $L_\infty$-homomorphisms corresponding to $\lambda$ and $\phi$.  The
$L_\infty$-homomorphism corresponding to $\lambda \maps \T{\Omega\g}
\to \Pg$ is given by
\begin{eqnarray*}
\lambda_0 (\ell) &=& \ell \cr
\lambda_1 (\ell) &=& (\ell, 0) \cr
\lambda_2 (\ell_1, \ell_2) &=&
\displaystyle{\big( 0,\, 2k \int_0^{2\pi}
\bracket{\ell_1}{\ell_2'} \, d\theta \big)}
\end{eqnarray*}
where $\ell, \ell_1, \ell_2 \in \Omega\g$, and that corresponding
to $\phi \maps \Pg \to \g_k$ is given by:
\begin{eqnarray*}
\phi_0 (p) &=& p(2\pi) \cr
\phi_1 (\ell, c) &=& c \cr
\phi_2 (p_1, p_2) &=& \displaystyle{k
       \int_0^{2\pi} \big( \bracket{p_2}{p_1'} -
                   \bracket{p_2'}{p_1} \big) \; d\theta }
\end{eqnarray*} where $p, p_1, p_2 \in P_0 \g$,
$\ell \in \Og$, and $c \in \R$.

\begin{lem} \et The composite
\[
    \T{\Omega\g}
    \stackto{\lambda}
    \Pg
    \stackto{\phi}
    \g_k
\]
is the zero homomorphism, and the kernel of $\phi$ is precisely the image
of $\lambda$, both on objects and on morphisms.
\end{lem}

\Proof
The composites $(\phi \circ \lambda)_0$ and
$(\phi \circ \lambda)_1$ clearly vanish.  Moreover
$(\phi \circ \lambda)_2$ vanishes since:
\begin{eqnarray*}
    (\phi \circ \lambda)_2 (\ell_1,\ell_2)
    &=&
    \phi_2 (\lambda_0 (\ell_1) , \lambda_0 (\ell_2) )
    +
    \phi_1 ( \lambda_2 ( \ell_1,\ell_2 )) \qquad
    \textrm{by} \; (\ref{composite.homo})
    \\
    &=&
    \phi_2 ( \ell_1,\ell_2 )
    +
    \phi_1 \big(0, 2k \int \bracket{\ell_1}{\ell_2'} \, d\theta \big)
    \\
    &=&
    k \int ( \bracket{\ell_2}{\ell_1'}
    - \bracket{\ell_2 '}{\ell_1} ) \, d\theta
    + 2k \int \bracket{\ell_1}{\ell_2'} \, d\theta  \\
    &=& 0
\end{eqnarray*}
with the help of integration by parts.
The image of $\lambda$ is precisely the kernel of $\phi$ by construction.
\endofproof

Note that $\phi$ is obviously onto, both for objects and
morphisms, so we have an exact sequence
\[    0 \to \T\Omega\g \stackto{\lambda}
            \Pg \stackto{\phi} \g_k  \to 0 \, . \]
Next we construct a family of splittings
$\psi \maps \g_k \to \Pg$ for this exact sequence:

\begin{lem} \et \label{psi}
Suppose
\[
    f \maps [0,2\pi] \to \mathbb{R}
\]
is a smooth function with $f(0) = 0$ and $f(2\pi) = 1$.
Then there is a Lie 2-algebra homomorphism
\[
   \psi \maps \g_k \to \Pg
\]
whose corresponding $L_{\infty}$-homomorphism is given by:
\begin{eqnarray*}
\psi_0 (x) &=& xf \cr
\psi_1 (c) &=& (0, c) \cr
\psi_2 (x_1, x_2) &=& ([x_1, x_2](f-f^2), 0)
\end{eqnarray*}
where $x, x_1, x_2 \in \g$ and $c \in \R$.
\end{lem}

\Proof
We show that $\psi$ satisfies the conditions of Definition \ref{Linftyhomo}.
We begin by showing that $\psi$ is a chain map, meaning that the following
square commutes:
$$\xymatrix{
    \R \ar[rr]^{d}
    \ar[dd]_{\psi_{1}} &&
    \g  \ar[dd]^{\psi_{0}} \\ \\
    \wOkg \ar[rr]^{d'} &&
    P_0 \g
}$$
where $d =0 $ since $\g_k$ is skeletal and $d'$ is the composite given
in Proposition \ref{Pg}.  This square commutes because
$\psi_0(d(c)) = \psi_0 (0) = 0$ and $d'(\psi_1(c)) =
d'(0,c) = 0$.

We continue by verifying conditions (\ref{homo1}) - (\ref{homo3})
of Definition \ref{Linftyhomo}.
The right-hand side of (\ref{homo1}) equals:
$$\psi_0 ( l_2(x_1,x_2) ) - l_2 ( \psi_0 (x_1) , \psi_0 (x_2) )
=  [x_1,x_2] (f-f^2) \, . $$
This equals the left-hand side $d'(\psi_2 (x_1,x_2)) $ by
construction.

The right-hand side of (\ref{homo2}) equals:
\[
  \psi_1( l_2(x,c) ) - l_2( \psi_0(x), \psi_1(c) ) \; =  \;
  \psi_1( 0 ) - l_2( xf, (0,c) ) \; = \; 0
\]
since both terms vanish separately. Since the left-hand side is
$\psi_2(x,dc) = \psi_2(x,0) = 0$, this shows that $\psi$ satisfies
condition (\ref{homo2}).

Finally we verify condition (\ref{homo3}).  The term
$l_3 ( \psi_0 (\cdot), \psi_0 (\cdot), \psi_0 (\cdot))$ vanishes
because $\Pg$ is strict. The sum of three other terms vanishes
thanks to the Jacobi identity in $\g$:
\begin{eqnarray*}
&& \psi_2(x_1, l_2(x_2,x_3))
 + \psi_2(x_2, l_2(x_3,x_1))
 + \psi_2(x_3, l_2(x_1,x_2))  \\
    &=&
    \Big(([x_1,[x_2,x_3]] +
          [x_2,[x_3,x_1]] +
          [x_3,[x_1,x_2]]) \, (f - f^2),\, 0 \Big)  \\
    &=&
    (0,0) \, .
\end{eqnarray*}
Thus, it remains to show that:
$$
 \psi_1 (l_3(x_1,x_2,x_3)) =  $$
$$  l_2(\psi_0 (x_1), \psi_2 (x_2, x_3))
+ l_2(\psi_0 (x_2), \psi_2 (x_3, x_1))
+ l_2(\psi_0 (x_3), \psi_2 (x_1, x_2)) \, .
$$
This goes as follows:
\begin{eqnarray*}
    &&   l_2(\psi_0 (x_1), \psi_2 (x_2, x_3))
       + l_2(\psi_0 (x_2), \psi_2 (x_3, x_1))
       + l_2(\psi_0 (x_3), \psi_2 (x_1, x_2)) \\
    &=&
    \Big(0, -3 \cdot 2k \int_0^{2\pi}
      \bracket{x_1}{[x_2,x_3]} \, f(f-f^2)' \, d\theta \; \Big) \\
    &=&
    \left(0,k\bracket{x_1}{[x_2,x_3]} \right) \qquad
\textrm{by the calculation below} \\
    &=&
    \psi_1 (l_3(x_1, x_2, x_3)) \, .
\end{eqnarray*}
The value of the integral here is \emph{universal},
independent of the choice of $f$:
\begin{eqnarray*}
\int_0^{2\pi}  f (f-f^2)' \; d \theta
&=&
\int_0^{2\pi} \left(f(\theta) f'(\theta)
- 2 f^2 (\theta) f'(\theta) \right) \; d\theta
    \\
    &=&
    \frac{1}{2} - \frac{2}{3}
    =
    -\frac{1}{6}
    \,.
\end{eqnarray*}
\hbox{\hskip 50ex} \endofproof

The final step in proving Theorem \ref{bigthm}
is to show that $\phi \circ \psi$ is the
identity on $\g_k$, while $\psi \circ \phi$ is isomorphic to the
identity on $\Pg$.
For convenience, we recall the definitions first:
$\phi \maps \Pg \to \g_k$ is given by:
\begin{eqnarray*}
\phi_0 (p) &=& p(2\pi) \cr
\phi_1 (\ell, c) &=& c \cr
\phi_2 (p_1, p_2) &=& \displaystyle{k
     \int_0^{2\pi}
\left( \bracket{p_2}{p_1'} - \bracket{p_2'}{p_1} \right) \; d\theta }
\end{eqnarray*}
where $p, p_1, p_2 \in P_0 \g$, $\ell \in \wOkg$, and $c \in \R$,
while $\psi \maps \g_k \to \Pg$ is given by:
\begin{eqnarray*}
\psi_0 (x) &=& xf \cr
\psi_1 (c) &=& (0, c) \cr
\psi_2 (x_1, x_2) &=& ([x_1, x_2] (f-f^2), 0)
\end{eqnarray*}
where $x, x_1, x_2 \in \g$, $c \in \R$, and
$f \maps [0,2\pi] \to \R$ satisfies the conditions of Lemma \ref{psi}.

\begin{lem} \et
      With the above definitions we have:
    \begin{itemize}
      \item $\phi \circ \psi$ is the identity Lie 2-algebra homomorphism
            on $\g_k$;
      \item $\psi \circ \phi$ is isomorphic, as a Lie 2-algebra
            homomorphism, to the identity on $\Pg$.
    \end{itemize}
\end{lem}

\Proof
We begin by demonstrating that $\phi \circ \psi$
is the identity on $\g_k$.
First,
$$(\phi \circ \psi)_0 (x) = \phi_0 (\psi_0 (x))
= \phi_0(xf) = xf(2\pi) = x,$$
since $f(2\pi) = 1$ by the definition of $f$ in Lemma \ref{psi}.
Second,
$$(\phi \circ \psi)_1 (c) = \phi_1 (\psi_1 (c))
= \phi_1 ((0,c)) = c$$
Finally,
\begin{eqnarray*}
(\phi \circ \psi)_2 (x_1, x_2)
&=&
\phi_2 (\psi_0 (x_1), \psi_0(x_2)) + \phi_1 (\psi_2 (x_1, x_2))
\qquad \textrm{by} \; (\ref{composite.homo}) \\
&=&
\phi_2(x_1 f, x_2 f) \; + \; \phi_1([x_1, x_2](f-f^2),0) \\
&=&
k \int (
\bracket{x_2 f}{x_1 f'} - \bracket{x_2 f'}{x_1 f} ) \, d\theta \;
+ \; 0 \\
&=&
k \bracket{x_2}{x_1} \int (ff' - f'f)\; d\theta
\\
&=& 0  \, .
\end{eqnarray*}

Next we consider the composite
\[
    \psi \circ \phi \maps \Pg \to \Pg \, .
\]
The corresponding $L_{\infty}$-algebra homomorphism is given by:
\begin{eqnarray*}
(\psi \circ \phi)_0 (p) &=& p(2\pi)f \\
(\psi \circ \phi)_1 (\ell,c) &=& (0,c) \\
(\psi \circ \phi)_2 (p_1,p_2) &=&
\left( [p_1 (2\pi) , p_2 (2\pi)](f-f^2), \;
k \int \left( \bracket{ p_2}{ p_1'} - \bracket{ p_2'}{ p_1}
\right) \; d\theta \right)
\end{eqnarray*}
where again we use equation (\ref{composite.homo}) to obtain
the formula for $(\psi \circ \phi)_2$.

For this to be isomorphic to the identity there must
exist a Lie 2-algebra 2-isomorphism
$$\tau \maps \psi \circ \phi \To \id$$
where $\id$ is the identity on $\Pg$.
We define this in terms of its corresponding
$L_{\infty}$-2-homomorphism by setting:
$$\tau(p) = (p-p(2\pi)f,0) \,.$$
Thus, $\tau$ turns a path $p$ into the loop $p-p(2\pi)f$.

We must show that $\tau$ is a chain homotopy satisfying
condition (\ref{2homo}) of Definition \ref{Linfty2homo}.
We begin by showing that $\tau$ is a chain homotopy.  We have
\begin{eqnarray*}
d(\tau(p)) \; = \; d(p-p(2\pi)f,0) &=& p -p(2\pi)f \\
&=& \id_0(p) - (\psi \circ \phi)_0 (p)
\end{eqnarray*}
and
\begin{eqnarray*}
\tau (d(\ell,c)) \; = \; \tau (\ell) &=& (\ell,0) \\
&=& \id_1(\ell, c) - (\psi \circ \phi)_1 (\ell, c)
\end{eqnarray*}
so $\tau$ is indeed a chain homotopy.

We conclude by showing that $\tau$ satisfies
condition (\ref{2homo}):
$$(\psi \circ \phi)_2 (p_1,p_2) =
l_2((\psi \circ \phi)_0 (p_1), \tau (p_2)) +
l_2 (\tau (p_1), p_2) - \tau ( l_2 (p_1, p_2))$$
In order to verify this equation, we write out the right-hand
side more explicitly by inserting the formulas for
$(\psi\circ\phi)_2$ and for $\tau$, obtaining:
\[
  l_2\big(p_1(2\pi)f,\, (p_2-p_2(2\pi)f,0)\big)
  +
  l_2\big((p_1-p_1(2\pi)f,0),\, p_2\big)
  -
  ([p_1,p_2] - [p_1(2\pi),p_2(2\pi)]f,\, 0)
\]
This is an ordered pair consisting of a loop in $\g$ and a real number.
By collecting summands, the loop itself turns out to be:
\[
  [p_1(2\pi),p_2(2\pi)] (f-f^2)
  \,.
\]
Similarly, after some integration by parts the real number is found to be:
\[
 k  \int_0^{2\pi}
    \left(
       \bracket{p_2}{p_1'}
       -
       \bracket{p_2'}{p_1}
    \right)
    \, d\theta
  \,.
\]
Comparing these results with the value of
$(\psi\circ \phi)_2(p_1,p_2)$ given
above, one sees that $\tau$ indeed satisfies (\ref{2homo}).
\endofproof


\section{Conclusions}
\label{conclusions.section}

We have seen that the Lie 2-algebra $\g_k$ is equivalent
to an infinite-dimensional Lie 2-algebra $\Pg$, and that when
$k$ is an integer, $\Pg$ comes from an infinite-dimensional
Lie 2-group $\PG$.  Just as the Lie 2-algebra $\g_k$ is built from
the simple Lie algebra $\g$ and a shifted version of $\u(1)$:
\[ 0 \stackto{\;} {\rm b}\u(1) \stackto{\;} \g_k \stackto{\;}
\g \stackto{\;} 0\, , \]
the Lie 2-group $\PG$ is built from $G$ and another Lie 2-group:
\[ 1 \stackto{\;} \LG \stackto{\;} \PG \stackto{\;} G \stackto{\;} 1  \]
whose geometric realization is a shifted version of $\U(1)$:
\[ 1 \stackto{\;} B\U(1) \stackto{\;}
|\PG| \stackto{\;} G \stackto{\;} 1\, . \] None of these exact
sequences split; in every case an interesting cocycle plays a role
in defining the middle term.  In the first case, the Jacobiator of
$\g_k$ is $k\nu \maps \Lambda^3 \g \to \R$. In the second case,
composition of morphisms is defined using multiplication in the
level-$k$ Kac--Moody central extension of $\OG$, which relies on
the Kac--Moody cocycle $k\omega \maps \Lambda^2 \Og \to \R$.  In
the third case, $|\PG|$ is the total space of a twisted
$B\U(1)$-bundle over $G$ whose Dixmier--Douady class is
$k[\nu/2\pi] \in H^3(G)$.  Of course, all these cocycles are
different manifestations of the fact that every simply-connected
compact simple Lie group has $H^3(G) = \Z$.

We conclude with some remarks of a more speculative nature.
There is a theory of `2-bundles' in which a Lie 2-group plays
the role of structure group \cite{BS,Bartels}.  Connections on
2-bundles describe parallel transport of 1-dimensional
extended objects, e.g.\ strings.  Given the importance of the
Kac--Moody extensions of loop groups in string theory, it is
natural to guess that connections on 2-bundles with structure
group $\PG$ will play a role in this theory.

The case when $G = \Spin(n)$ and $k = 1$ is particularly interesting,
since then $|\PG| = \String(n)$.  In this case we suspect that
$2$-bundles on a spin manifold $M$ with structure $2$-group $\PG$
can be thought as substitutes for principal $\String(n)$-bundles on
$M$.  It is interesting to think about `string structures'
\cite{MuSt} on $M$ from this perspective: given a principal $G$-bundle
$P$ on $M$ (thought of as a $2$-bundle with only identity morphisms)
one can consider the obstruction problem of trying to lift the structure
$2$-group from $G$ to $\PG$.  There should be a single topological
obstruction in $H^4(M;\Z)$ to finding a lift, namely the characteristic
class $p_1/2$.  When this characteristic class vanishes, every principal
$G$-bundle on $M$ should have a lift to a $2$-bundle $\mathcal{P}$ on $M$
with structure $2$-group $\PG$.  It is tempting to conjecture that the
geometry of these $2$-bundles is closely related to the enriched
elliptic objects of Stolz and Teichner \cite{Stolz-Teichner}.

\subsection*{Acknowledgements}
We thank Andr\'e Henriques for useful correspondence.
We thank Edward Witten for suggesting that we find a 2-group
related to elliptic cohomology.  We thank Aaron Lauda for the
beautiful picture of $\PG$.  We also thank David Roberts for catching
and correcting many sign mistakes in this, the version of
this paper placed on the arXiv in May 2023.


\begin{thebibliography}{10}

\bibitem{HDA5}
J.~Baez and A.~Lauda, Higher-dimensional algebra {V}: 2-groups,
{\sl Theory and Applications of Categories} {\bf 12} (2004), 423--491.
Also available as \href{https://arxiv.org/abs/math.QA/0307200}{arXiv:math.QA/0307200}.

\bibitem{HDA6}
J.~Baez and A.~Crans, Higher-dimensional algebra {VI}: Lie 2-algebras,
{\sl Theory and Applications of Categories} {\bf 12} (2004), 492--528.
Also available as  \href{https://arxiv.org/abs/math.QA/0307263}{arXiv:math.QA/0307263}.

\bibitem{BCSS}
J.~Baez, A.~Crans, D.~Stevenson and U.~Schreiber, From loop groups
to 2-groups, first version available as \href{https://arxiv.org/abs/math.QA/arXiv:0504123v1}{math.QA/0504123v1}.

\bibitem{BS}
J.~Baez and U.~Schreiber, Higher gauge theory, in
{\sl Categories in Algebra, Geometry and Mathematical
Physics}, eds. A. Davydov {\it et al},
Contemp.\ Math.\ {\bf 431}, AMS, Providence, Rhode Island, 2007.
Also available as \href{https://arxiv.org/abs/math.DG/0511710}{arXiv:math.DG/0511710}.

\bibitem{Bartels}
T.\ Bartels, Higher gauge theory I: 2-bundles, available as
\href{https://arxiv.org/abs/math.CT/0410328}{arXiv:math.CT/0410328}.

\bibitem{Brylinski} J.--L.\ Brylinski,
{\sl Loop Spaces, Characteristic Classes and Geometric Quantization},
Birkhauser, Boston, 1993.

\bibitem{Brylinski2} J.--L.\ Brylinski, Differentiable cohomology of
gauge groups, available as \href{https://arxiv.org/abs/math.DG/0011069}{arXiv:math.DG/0011069}.

\bibitem{BM} J.--L.\ Brylinski and D.\ A.\ McLaughlin,
The geometry of degree--four characteristic classes and of line
bundles on loop spaces I, {\sl Duke Math.\ J.\ }{\bf 75} (1994), 603--638.
II, {\sl Duke Math.\ J.\ }{\bf 83} (1996), 105--139.

\bibitem{CJMSW} A.\ Carey, S.\ Johnson, M.\ Murray, D.\ Stevenson
and B.--L.\ Wang, Bundle gerbes for Chern-Simons and Wess-Zumino-Witten
theories, to appear in {\sl Commun.\ Math.\ Phys.}  Also
available as \href{https://arxiv.org/abs/math.DG/0410013}{arXiv:math.DG/0410013}.

\bibitem{CheegerS} J.\ Cheeger and J.\ Simons, Differential
characters and geometric invariants, in {\sl Geometry and
Topology}, eds.\ J.\ Alexander and J.\ Harer,
Lecture Notes in Mathematics {\bf 1167} (1985), 50--80.

\bibitem{ChernS} S.\ S.\ Chern and J.\ Simons, Characteristic
forms and geometric invariants, {\sl Ann.\ Math.\ }{\bf 99} (1974),
48--69.

\bibitem{JS} A. \ Joyal and R. \ Street, Braided monoidal
categories, Macquarie Mathematics Report No. \ 860081, November
1986.

\bibitem{Hamilton} R.\ S.\ Hamilton, The inverse function theorem
of Nash and Moser, {\sl Bull.\ Amer.\ Math.\ Soc.\ }{\bf 7} (1982),
65--222.

\bibitem{Henriques} A.\ Henriques, Integrating $L_\infty$-algebras,
available as 
\href{https://arxiv.org/abs/math.AT/0603563}{arXiv:math.AT/} \break
\href{https://arxiv.org/abs/math.AT/0603563}{0603563}.

\bibitem{LM} T.\ Lada and M.\ Markl, Strongly homotopy Lie algebras,
{\sl Comm. Alg.\ }{\bf 6} (1995), 2147--2161. Also available
as \href{https://arxiv.org/abs/hep-th/9406095}{arXiv:hep-th/9406095}.

\bibitem{Mac} S.\ Mac Lane, {\sl Categories for the Working
Mathematician,} 2nd edition, Springer, Berlin, 1998.

\bibitem{Milnor} J.~Milnor, Remarks on infinite dimensional Lie groups,
{\sl Relativity, Groups and Topology II, (Les Houches, 1983)},
North-Holland, Amsterdam, 1984.

\bibitem{Mickelsson} J.\ Mickelsson,
Kac--Moody groups, topology of the Dirac determinant bundle,
and fermionization, {\sl Commun.\ Math.\ Phys.\ }{\bf 110}
(1987), 173--183.

\bibitem{Murray} M.\ K.\ Murray, Bundle gerbes,
{\sl J.\ Lond.\ Math.\ Soc.\ }{\bf 54} (1996), 403--416.
Also available as \href{https://arxiv.org/abs/dg-ga/9407015}{arXiv:dg-ga/9407015}.

\bibitem{Murray1} M.\ K.\ Murray, Another construction of the central 
extension of the loop group, {\sl Commun.\ Math.\ Phys.\ }{\bf 116} 
(1988), no.~1 73--80.

\bibitem{MuSt} M.\ K.\ Murray, D.\ Stevenson, Higgs fields, bundle
gerbes and string structures, {\sl Commun.\ Math.\ Phys.\ }{\bf 243}
(2003), 541--555.  Also available as \href{https://arxiv.org/abs/arXiv:math.DG/0106179}{math.DG/0106179}.

\bibitem{PressleySegal} A.~Pressley and G.~Segal, {\em Loop Groups}.
Oxford U.\ Press, Oxford, 1986.

\bibitem{SS} M.\ Schlessinger and J.\ Stasheff, The Lie algebra
structure of tangent cohomology and deformation theory, {\sl Jour.\
Pure App.\ Alg.\ }{\bf 38} (1985), 313--322.

\bibitem{Segal} G.\ B.\ Segal, Classifying spaces and spectral
sequences, {\sl Publ.\ Math.\ IHES} {\bf 34} (1968), 105--112.

\bibitem{Stacey} A.\ Stacey, personal communication.  

\bibitem{Stolz-Teichner} S.\ Stolz and P.\ Teichner,
What is an elliptic object?, in {\sl Topology, Geometry and
Quantum Field Theory: Proceedings of the 2002 Oxford Symposium
in Honour of the 60th Birthday of Graeme Segal}, ed.\ U.\ Tillmann,
Cambridge U.\ Press, Cambridge, 2004.

\bibitem{Witten} E.\ Witten, The index of the Dirac operator in
loop space, in {\sl Elliptic Curves and Modular Forms in Algebraic
Topology}, ed.\ P.\ S.\ Landweber, Lecture Notes in
Mathematics {\bf 1326}, Springer, Berlin, 1988, pp.\ 161--181.

\end{thebibliography}
\end{document}